\documentclass[11pt]{article}

\usepackage{bbm}
\usepackage{amsmath,amsfonts,amssymb,enumerate}
\usepackage{color}
\usepackage{lscape}
\usepackage{latexsym}
\usepackage{multirow}
\usepackage{rotating}
\usepackage{threeparttable}
\usepackage{graphicx}
\usepackage{algpseudocode,algorithm}
\usepackage[margin=1in]{geometry}

\numberwithin{equation}{section}
\newtheorem{theorem}{\bf Theorem}[section]
\newtheorem{lemma}[theorem]{\bf Lemma}
\newtheorem{proposition}[theorem]{\bf Proposition}

\newtheorem{coro}[theorem]{\bf Corollary}

\newenvironment{proof}{\noindent{\em Proof:}}{\quad \hfill$\Box$\vspace{2ex}}

\graphicspath{{figs/}}

\begin{document}

\title{\bf Sparse Regularization with the $\ell_0$ Norm\thanks{This work is supported in part by the US National Science Foundation under grant DMS-1912958}}

\author{Yuesheng Xu\thanks{Department of Mathematics and Statistics, Old Dominion University, Norfolk, Virginia 23529, USA, y1xu@odu.edu}}
\date{}

\maketitle

\begin{abstract}
We consider a minimization problem whose objective function is the sum of a fidelity term, not necessarily convex, and a regularization term defined by a positive regularization parameter $\lambda$ multiple of the $\ell_0$ norm composed with a linear transform. This problem has wide applications in compressed sensing, sparse machine learning and image reconstruction. The goal of this paper is to understand what choices of the regularization parameter can dictate the level of sparsity under the transform for a global minimizer of the resulting regularized objective function. This is a critical issue but it has been left unaddressed. We address it from a geometric viewpoint with which the sparsity partition of the image space of the transform is introduced. Choices of the regularization parameter are specified to ensure that a global minimizer of the corresponding regularized objective function achieves a prescribed level of sparsity under the transform. Results are obtained for the spacial sparsity case in which the transform is the identity map, a case that covers several applications of practical importance, including machine learning, image/signal processing and medical image reconstruction.
\end{abstract}



\section{Introduction}
The aim of this work is to understand a global minimizer of regularization problems whose objective functions have the form of a fidelity term plus a regularization term involving the $\ell_0$ norm. Regularization problems of this type appear frequently in recent studies of machine learning \cite{Lopez11, Louizos18, Pan, Wang2016}, computer graphics  \cite{He13, Sun15}, signal processing \cite{Dai16, Liu18, Zeng19}, image processing \cite{Shen16, Shen14, Zeng18}, medical imaging \cite{Zheng19} and statistics \cite{Huang18, Zheng14}. Many published results have demonstrated that the use of the $\ell_0$ norm in regularization models promotes sparsity for the regularized solutions or the transformed regularized solutions. Most of the existing work focus on developing numerical algorithms and considering convergence issues of the developed algorithms. It remains to be understood how choices of the regularization parameter balance the sparsity of a global minimizer of the regularization problem and its approximation to a global minimizer of the fidelity function. It is the goal of this paper to provide mathematical understanding on how the use of the $\ell_0$ norm as a regularization term promotes sparsity of the regularized solutions or the transformed regularized solutions.

We now describe precisely the problem to be considered in this paper.
Let $d$ be a fixed positive integer. For $x\in \mathbb{R}^d$, we use $\|x\|_0$
to denote the number of the nonzero components of $x$. Although $\|\cdot\|_0$ does not satisfy the axiom of vector norms, it is widely referred to as the $\ell_0$ norm in the sparse optimization community. We follow the custom of the community to call it the $\ell_0$ norm. Let $m$ be another positive integer, which may be equal to $d$ or may be different from $d$. Suppose that $g: \mathbb{R}^m\to \mathbb{R}$ is a given function and $M$ is a real $d\times m$ matrix. For a parameter $\lambda>0$, we define the function
\begin{equation}\label{main-function-transform}
f(x):=g(x)+\lambda\|Mx\|_0, \ \ x\in\mathbb{R}^m
\end{equation}
and consider the related regularization problem
\begin{equation}\label{regularization}
\min \{f(x):\ x\in\mathbb{R}^m\}.
\end{equation}
Here, $\lambda>0$ is a regularization parameter. Its choices may impose sparsity of a global minimizer of the corresponding function \eqref{main-function-transform}.
Clearly, the function $f$ defined by equation \eqref{main-function-transform} depends on the parameter $\lambda$ and the transform matrix $M$. Although for  conciseness of notation, we do not label the dependence of $f$ on $\lambda$ or $M$ in its notation, we always assume that $f$ depends on these quantities.

In the context of regularization, the function $g$ appearing in \eqref{main-function-transform} is the data fidelity term derived from a linear \cite{Kress} or nonlinear ill-posed problem \cite{Hofmann}. It may also describes a network \cite{Le} in machine learning. For more linear ill-posed problems, see \cite{JChen, CMX, ChenXuYang}. The function $g$ that appears in application is often convex, (for example, the least squares error). It can also be non-convex. For instance, fidelity terms for deep learning are non-convex \cite{Unser}. It can be differentiable or non-differentiable. In this paper, in order to enlarge the applicability of the established theory, we consider a wide class of fidelity terms $g$, without imposing convexity or differentiability.

The matrix $M$ that appears in the regularization term is often chosen as a mathematical transform such as a discrete cosine transform \cite{Strang}, a wavelet transform \cite{Daub, Lian11, Mallat, MX94, MX97} or a framelet transform \cite{Chan, Ron}, depending on specific applications. It can also be a difference matrix (for example, the $\ell_0$-TV). For TV-regularization, the readers are referred to \cite{Osher}. The matrix $M$ does not have to be a square matrix. However, we confine ourselves to matrices of full rank, since most of mathematical transforms used frequently in applications have this property and the case with matrices of arbitrary rank may be treated by employing the singular value decomposition, on which we will comment at the end of the last section.

The regularization problem \eqref{regularization} often raises in the scenarios that the function $g$ has a global minimizer which itself may not be sparse while a sparse minimizer is desirable. Bringing forward such a model enables us to find a global minimizer of $f$ having the desired sparsity under the transform while keeping it as close to the global minimizer of $g$ as possible. A desirable solution of the problem \eqref{regularization} is the one achieving the desired sparsity and being close to the global minimizer of the fidelity term $g$. For this reason, we shall assume that the function $g$ has a global minimizer in $\mathbb{R}^m$.

We are also interested in minimization problems of sparsity regularization in the spacial domain, that is, the special case of \eqref{main-function-transform} with $d=m$ and $M=I$, the identity matrix. In this case, the function $f$ has the spacial form
\begin{equation}\label{main-function}
f(x):=g(x)+\lambda\|x\|_0, \ \ x\in \mathbb{R}^d.
\end{equation}
Although the model \eqref{main-function} has its practical importance, we shall not present special results for this case since they can be obtained from general results by restricting $M=I$.

Motivated from {\it approximately sparse} regularization such as regularization by the envelope of the $\ell_0$-norm and capped-$\ell_1$, we introduce the function
\begin{equation}\label{def:f}
f(x,y):=g(x,y)
+\lambda\|x\|_0,\ \ \mbox{for all}\ \ (x, y)\in \mathbb{R}^d\times \mathbb{R}^{d'},
\end{equation}
where $g: \mathbb{R}^d\times \mathbb{R}^{d'}\to \mathbb{R}$ is a continuous function, not necessarily convex. Typical examples of function $f$ in the form \eqref{def:f} include the objective functions in wavelet inpainting with the $\ell_0$ sparse regularization \cite{Shen16}, inverting an incomplete Fourier transform \cite{WuXu2021} and medical image reconstruction \cite{Zheng19}. In these applications, the function $g$ is a sum of two or three convex functions which measure the data fidelity and define other convex constraints. We shall study what choices of the positive parameter $\lambda$ will balance the sparsity of global minimizers of function $f$ and its approximation to global minimizers of function $g$. We do not intend to provide practical methods for choices of the regularization parameter $\lambda$, and rather, we supply a mathematical understanding of the relation among choices of the regularization parameter, and global/local minimizers of the two functions $f$ and $g$.  We are also interested in understanding the relation between local minimizers of these two functions when the regularization parameter $\lambda$ is fixed.

Many existing empirical results exhibit that when the regularization parameter $\lambda$ is appropriately chosen, a global minimizer of a function $f$ in the form \eqref{main-function-transform} has certain sparsity under the given transform $M$ and a global minimizer of $f$ in the form \eqref{main-function} has sparsity in the spacial domain. We are interested in understanding mathematically how the sparsity of a global minimizer of the function $f$ having the form \eqref{main-function-transform} or  \eqref{main-function} depends on choices of the parameter $\lambda$.
For convenience of presentation, we shall refer the case with a general matrix $M$ as to sparsity under the transform and the special case when $M=I$ to sparsity in the spacial domain, when clarification is desired.

Our key approach is the understanding of the ``surface'' geometry of the function $f$ defined by \eqref{main-function-transform} or \eqref{main-function}. When $\lambda=0$, $f$ clearly reduces to the function $g$. We regard the surface determined by the function $g$ as the original ``landscape" and imagine its animation controlled by the parameter $\lambda>0$. At the moment when we start to increase the value of the parameter $\lambda$ from $0$ to a positive number, the original ``landscape" begins to change like vertical fractures of the earth crust during an earthquake. The parts of the landscape corresponding to $Mx=0$ will stay in their original positions and other parts will lift upward according to $\lambda\|Mx\|_0$. This geometry motivates us to partition the space $\mathbb{R}^d$ according to the values of the $\ell_0$ norm of the vectors in the space, that is, the sparsity levels. This sparsity partition of the Euclidean space will enable us to understand how the value of the parameter $\lambda$ will determine the sparsity level of a global minimizer of $f$. We shall introduce the sparsity partition of the space $\mathbb{R}^d$, the image space of the transform $M$, and understand how this partition will result in a partition of the preimage space $\mathbb{R}^m$. Through these partitions we shall be able to visualize the animation as the value of the parameter $\lambda$ increases. As a result, we can clearly determine how large the value of $\lambda$ will be in order to achieve a desired level of sparsity for a global minimizer of $f$ and at the same time to keep the minimizer as close to the global minimizer of the function $g$ as possible.

We organize this paper in five sections. In section 2, we introduce a partition of the image space of a transform $M$ according to the levels of sparsity and consider its corresponding partition of the preimage space. We study both algebraic and topological properties of the sets in these partitions. We devote section 3 to a study of choices of the parameter $\lambda>0$ that ensure desired levels of sparsity under the transform of a global minimizer of function $f$ having the form \eqref{main-function-transform}. Several necessary conditions of a global minimizer of $f$ are presented. In section 4, for functions $f$ having the form \eqref{def:f}, we investigate the same issues as those considered in section 3 for function $f$ in the form \eqref{main-function-transform}.
We also present a relation between local minimizers of minimization problem \eqref{opt:f} and its reduced minimization problem without the term involving the $\ell_0$-norm. We briefly discuss in section 5 extension of the results presented in  section 3 and  section 4 involving matrix $M$ and potential practical uses of the main results of this paper.

\section{Sparsity Partition of the Euclidean Space}

We introduce in this section a partition of the space $\mathbb{R}^d$, the image space of the linear transform $M$, according to levels of sparsity, and study its corresponding partition of the preimage space $\mathbb{R}^m$. For the purpose of understanding the sparsity of a global minimizer of the function $f$ defined by \eqref{main-function-transform}, we present algebraic and topological properties of the sets in the partitions.

It is convenient to introduce the level of sparsity for a vector in $\mathbb{R}^d$. To this end, for a positive integer $d$, we define two index sets
$\mathbb{N}_d:=\{1,2, \dots, d\}$ and $\mathbb{Z}_{d}:=\{0, 1, \dots, d-1\}$.
Precisely, a vector $x\in \mathbb{R}^d$ is said to have sparsity
of level $\ell\in \mathbb{Z}_{d+1}$ if $x$ has exactly $\ell$ number of
nonzero components. Clearly, the zero vector has sparsity of level $0$ and a vector whose components are all nonzero have sparsity of level $d$. Vectors having sparsity of level $d$ are not sparse. Sparse vectors are those located on the coordinate axes or coordinate planes of space $\mathbb{R}^d$. For example, in $\mathbb{R}^3$, vectors on the three coordinate axes but not at the origin have sparsity of level 1, vectors on the three coordinate planes but not on the three coordinate axes have sparsity of level 2 and vectors not on the three coordinate planes have sparsity of level 3. Most vectors in the space $\mathbb{R}^d$ are not sparse. In fact, the set of the sparse vectors in $\mathbb{R}^d$ has zero measure.


We now define the sparsity partition of $\mathbb{R}^d$. We need the canonical basis for the space $\mathbb{R}^d$. For each $j\in\mathbb{N}_d$, by $e_j\in \mathbb{R}^d$, we
denote the unit vector with 1 for the $j$-th component and 0 otherwise. The vectors $e_j$, $j\in \mathbb{N}_d$, form the canonical basis for $\mathbb{R}^d$. Let

\begin{align}\label{setA}
& A_0:= \{0\in\mathbb{R}^d\}, \\
& A_\ell:=\left\{\sum_{j\in \mathbb{N}_\ell}x_{k_j}e_{k_j}: x_{k_j}\in
\mathbb{R}\setminus \{0\}, \ \mbox{for}\ 1\leq k_1<k_2<\cdots<k_\ell\leq
d \right\},  \ \ \mbox{for}\ \ \ell\in
\mathbb{N}_d.\nonumber
\end{align}

In the next proposition, we show that the sets $A_\ell$, $\ell\in \mathbb{Z}_{d+1}$, defined by \eqref{setA} indeed form a partition for the space $\mathbb{R}^d$.

\begin{proposition}\label{sets}
If the sets $A_\ell$, $\ell\in \mathbb{Z}_{d+1}$, are defined by \eqref{setA}, then

(i) they  are mutually disjoint,

(ii) they form a partition for the space $\mathbb{R}^d$, that is,
\begin{equation}\label{decompositionOfRd}
\mathbb{R}^d=\bigcup_{\ell\in \mathbb{Z}_{d+1}}A_\ell.
\end{equation}
\end{proposition}
\begin{proof}
(i) It suffices to show that
\begin{equation}\label{intercept}
A_j\cap A_{j'}=\emptyset, \ \ \mbox{for all}\ \ j,j'\in \mathbb{Z}_{d+1}\ \ \mbox{with}\ \ j\neq j'.
\end{equation}
Without loss of generality, we assume that $j<j'$. Suppose that $x\in A_j\cap A_{j'}$. By the definition of $A_j$, there exist $1\leq k_1<k_2<\cdots<k_j\leq d$ and $x_{k_i}\in\mathbb{R}\setminus \{0\}$ such that
\begin{equation}\label{x-j}
x=\sum_{i\in \mathbb{N}_j}x_{k_i}e_{k_i},
\end{equation}
and by the definition of $A_{j'}$, there exist $1\leq k'_1<k'_2<\cdots<k'_{j'}\leq d$ and $x'_{k'_i}\in
\mathbb{R}\setminus \{0\}$ such that
\begin{equation}\label{x-j'}
x=\sum_{i\in \mathbb{N}_{j'}}x'_{k'_i}e_{k'_i}.
\end{equation}
Subtracting equation \eqref{x-j} from \eqref{x-j'} yields
\begin{equation}\label{x-x-j'}
\sum_{i\in \mathbb{N}_{j'}}x'_{k'_i}e_{k'_i}-\sum_{i\in \mathbb{N}_j}x_{k_i}e_{k_i}=0
\end{equation}
We introduce two index sets
$\mathbb{I}:=\{k_1, k_2, \dots, k_j\}$ and  $\mathbb{I}':=\{k'_1, k'_2, \dots, k'_{j'}\}$.
Since $j<j'$, we observe that  $\mathbb{I}\neq  \mathbb{I}'$.
It follows that there exists an index $k'_t\in \mathbb{I}'$ but  $k'_t\notin \mathbb{I}$.
Since $e_j$, $j\in \mathbb{N}_d$, are linearly independent, according to \eqref{x-x-j'}, we conclude that $x'_{k'_t}=0$. This contradicts the hypothesis that $x'_{k'_t}\neq 0$ and confirms \eqref{intercept}.

(ii) Assume that $x\in\mathbb{R}^d$. Let $\ell:=\|x\|_0$. Then, $\ell\in \mathbb{Z}_{d+1}$. Thus, we have that $x\in A_\ell$. This ensures that  
$$
\mathbb{R}^d\subseteq\bigcup_{\ell\in \mathbb{Z}_{d+1}}A_\ell.
$$
Clearly, we have that 
$$
\bigcup_{\ell\in \mathbb{Z}_{d+1}}A_\ell\subseteq \mathbb{R}^d.
$$
These two inclusions imply the validity of equation \eqref{decompositionOfRd}, which together with part (i) of this proposition confirms that the sets $A_\ell$, $\ell\in \mathbb{Z}_{d+1}$, form a partition for $\mathbb{R}^d$.
\end{proof}

We illustrate Proposition \ref{sets} by $\mathbb{R}^2$. Clearly, for $d=2$, $\mathbb{R}^2=A_0\cup A_1\cup A_2$, where
$$
A_0:=\{(0,0)\}, \ \ A_1:=\{(x,0): x\in\mathbb{R}\setminus \{0\}\}\cup\{(0,y): y\in \mathbb{R}\setminus\{0\}\},
$$
and
$$
A_2:=\{(x,y): (x,y)\in \mathbb{R}\times\mathbb{R},\ x\neq 0 \ \mbox{and} \ y\neq 0\}.
$$
That is, $A_1$ contains points on the two axes except the origin and $A_2$ contains the four quadrants of the two dimensional plane.

We remark that according to \eqref{setA}, for each $\ell\in \mathbb{Z}_{d+1}$, $A_\ell$ is the set of all vectors in $\mathbb{R}^d$ having sparsity of level $\ell$.
According to Proposition \ref{sets}, the space $\mathbb{R}^d$ has the sparsity partition $A_j$, $j\in \mathbb{Z}_{d+1}$, which groups the vectors in $\mathbb{R}^d$ according to their sparsity levels.
We further observe that the sets $A_\ell$, $\ell\in \mathbb{N}_d$, are closed under the operation of nonzero scalar multiplication, but not closed under the operation of addition. For example, $e_1, e_2\in A_1$ but $e_1+e_2\in A_2$.

It is also convenient to define the set of vectors in $\mathbb{R}^d$ whose sparsity levels do not exceed $\ell$.
For $\ell\in \mathbb{Z}_{d+1}$, we let
$$
\Omega_\ell:=\bigcup_{j\in\mathbb{Z}_{\ell+1}}A_j.
$$
Clearly, $\Omega_\ell$ is the set of vectors in $\mathbb{R}^d$ whose sparsity levels do not exceed $\ell$. Moreover, we have that
\begin{equation}\label{Omega}
\Omega_0=A_0,
\ \
\Omega_{j+1}=\Omega_j\cup A_{j+1}, \ \ j\in \mathbb{Z}_{d+1}
\ \
\mbox{and}\ \
\Omega_d=\mathbb{R}^d.
\end{equation}
These equations yield that
$$
A_d=\mathbb{R}^d\setminus \Omega_{d-1}.
$$
The set $A_d$ consists of the vectors in $\mathbb{R}^d$ whose components are all nonzero.
By the definition of the sets $\Omega_j$ and properties of $A_j$, we see that $\Omega_j$ for $j\in \mathbb{N}_{d-1}$ are closed under the operation of nonzero scalar multiplication, but not closed under the operation of addition.

We now consider a partition of the space $\mathbb{R}^m$, the preimage space of the transform $M$, induced by the sparsity partition of $\mathbb{R}^d$.
Suppose that 
\begin{equation}\label{MR}
M\mathbb{R}^m=\mathbb{R}^d.
\end{equation}
When condition \eqref{MR} is satisfied, we say that $M$ is of full rank.
We introduce $d+1$ subsets $B_j$, $j\in \mathbb{Z}_{d+1}$, of the preimage space $\mathbb{R}^m$ according to the sparsity partition $A_j$, $j\in \mathbb{Z}_{d+1}$ by
$$
B_j:=\{x\in \mathbb{R}^m: Mx\in A_j\}.
$$
Because $A_0=\{0\}$, the set $B_0$ is the null space of matrix $M$. Moreover, we have the following simple fact.

\begin{proposition}\label{A-B}
If $M$ is a $d\times m$ full rank matrix, then
$$
MB_j=A_j, \ \ \mbox{for all}\ \ j\in \mathbb{Z}_{d+1}.
$$
\end{proposition}
\begin{proof}
Let $j\in \mathbb{Z}_{d+1}$ be fixed. We assume that $y\in MB_j$. Thus, there exists $x\in B_j$ such that $y=Mx$. By the definition of $B_j$, we have that $Mx\in A_j$. Hence, $y\in A_j$. This implies the inclusion $MB_j\subseteq A_j$. 

Conversely, we let $y\in A_j$. By Proposition \ref{sets}, the sets $A_j$, $j\in \mathbb{R}^d$, form a partition for the space $\mathbb{R}^d$ and thus, $y\in \mathbb{R}^d$. Since $M$ is of full rank, according to equation \eqref{MR}, there exists $x\in \mathbb{R}^m$ such that $y=Mx$. Since $Mx\in A_j$, we find that $x\in B_j$. Thus, we have that $y\in MB_j$. This yields the inclusion $A_j\subseteq MB_j$. We therefore establish the desired equation of this proposition.
\end{proof}

Proposition \ref{A-B} clearly reveals that for each $j\in \mathbb{Z}_{d+1}$, the set $B_j$ is the preimage set of $A_j$, the set of the vectors in $\mathbb{R}^d$ having sparsity of level $j$, under the transform $M$. However, vectors in $B_j$ do not necessarily have sparsity of level $j$.
In the next proposition, we show that the sets $B_j$, $j\in \mathbb{Z}_{d+1}$, form a partition for the preimage space $\mathbb{R}^m$ of the transform $M$.

\begin{proposition}
If $M$ is a $d\times m$ full rank matrix, then the sets $B_j$, $j\in \mathbb{Z}_{d+1}$, form a partition for the space $\mathbb{R}^m$.
\end{proposition}
\begin{proof}
It suffices to establish that
\begin{equation}\label{partitionB}
\mathbb{R}^m=\bigcup_{j\in\mathbb{Z}_{d+1}}B_j
\end{equation}
and
\begin{equation}\label{B-inters}
B_j\cap B_{j'}=\emptyset \ \ \mbox{for all}\ \ j,j'\in \mathbb{Z}_{d+1}\ \ \mbox{with}\ \ j\neq j'.
\end{equation}
To show \eqref{partitionB}, we let $x\in \mathbb{R}^m$. By the hypothesis on matrix $M$, we see that equation \eqref{MR} holds and thus, $Mx\in \mathbb{R}^d$. Employing the sparsity partition $A_j$, $j\in \mathbb{Z}_{d+1}$, for the space $\mathbb{R}^d$, we see that there exists $j\in \mathbb{Z}_{d+1}$ such that $Mx\in A_j$. By the definition of the set $B_j$, we conclude that $x\in B_j$. Hence, we have that
$$
\mathbb{R}^m\subseteq\bigcup_{j\in\mathbb{Z}_{d+1}}B_j.
$$
By the definition of the sets $B_j$, each of these sets is contained in $\mathbb{R}^m$. Thus, 
$$
\bigcup_{j\in\mathbb{Z}_{d+1}}B_j\subseteq\mathbb{R}^m.
$$
Consequently,
equation \eqref{partitionB} holds true.

It remains to prove equation \eqref{B-inters}. Suppose that $x\in B_j\cap B_{j'}$ for a fixed pair of indices $j,j'\in\mathbb{Z}_{d+1}$, with $j\neq j'$.  By the definition of the set $B_j$, we have that $Mx\in A_j$ and by the definition of the set $B_{j'}$, we have that $Mx\in A_{j'}$. According to Proposition \ref{sets}, the two sets $A_j$ and $A_{j'}$ are disjoint. This clearly implies that $Mx\in\emptyset$. Noting that $M$ is of full rank, we conclude that $x\in\emptyset$. That is, equation \eqref{B-inters} holds true.
\end{proof}

In the remaining part of this section, we study useful topological properties of the sets that we introduced earlier in this section.

\begin{proposition}\label{Closedness}
The following statements holds true.

(i) The set $A_0$ is closed, the sets $A_\ell$, $\ell\in \mathbb{N}_{d-1}$, are neither closed nor open, and $A_d$ is open.

(ii) For $j\in \mathbb{Z}_d$, $\Omega_j$ are closed sets.
\end{proposition}
\begin{proof}
(i) Since $A_0$ contains only one point $0$, it is closed. It is straightforward to see that the sets $A_j$, $j\in \mathbb{N}_{d-1}$, are not open since in every neighbourhood of a vector in $A_j$ contains vectors that are not in $A_j$. We now show that the sets $A_j$, $j\in \mathbb{N}_{d-1}$, are not closed either. To this end, we consider a sequence of vectors $x_n$, $n\in \mathbb{Z}$, in $\mathbb{R}^d$ whose first $j-1$ components are all equal to 1, last $d-j$ components are all equal to zero and $j$th component is $1/n$, that is,
$$
x_n:=(1, \dots, 1, 1/n, 0, \dots, 0), \ \ \mbox{for all}\ \ n\in \mathbb{Z}.
$$
Clearly, $x_n\in A_j$, for all $n\in \mathbb{Z}$ and $\|x_n-\hat x\|_2=1/n\to 0$ as $n\to \infty$, where $\hat x\in \mathbb{R}^d$ whose first $j-1$ components are all equal to 1 and last $d-j+1$ components are all equal to zero, that is, $\hat x:= (1, \dots, 1, 0, 0, \dots, 0)$. In other words, $x_n$ converges to a vector in $A_{j-1}$ not in $A_j$. Therefore, $A_j$, $j\in \mathbb{N}_{d-1}$, are not closed.

It remains to show that $A_d$ is open. Suppose that $\hat{x}\in A_d$. Then, we have that $\hat{x}=(\hat{t}_1, \hat{t}_2, \dots, \hat{t}_d)$ with $\hat{t}_j\neq 0$, for all $j\in\mathbb{N}_d$. Hence, for all  $j\in\mathbb{N}_d$ there exists $\epsilon>0$ such that for all $t'_j\in (\hat{t}_j-\epsilon/d^{1/2}, \hat{t}_j+\epsilon/d^{1/2})$, we have that $t'_j\neq 0$. Let $x':=(t'_1, t'_2, \dots, t'_d)$. We observe that $x'\in A_d$. That is, the open ball
$$
B_o(\hat{x}, \epsilon):=\{x\in \mathbb{R}^d: \|x-\hat{x}\|_2<\epsilon\}
$$
is contained in $A_d$. Thus, $A_d$ is an open set.

(ii) For a fixed $j\in \mathbb{Z}_d$, we assume that a sequence $x_n$, $n\in \mathbb{Z}$, in $\Omega_j$ converges to a point $x\in \mathbb{R}^d$ as $n\to \infty$, and we show that $x\in \Omega_j$ by contradiction. Assume, to the contrary, that $x\notin \Omega_j$. By the second equation of \eqref{Omega}, we have that $\Omega_j\subset\Omega_{j+1}$. Without loss of generality, we assume that $x\in \Omega_{j+1}$. Hence, $x\in A_{j+1}$. That is, $x$ has exactly $j+1$ nonzero components. Therefore, for sufficiently large $n$, $x_n$ has at least $j+1$ nonzero components. This contradicts the assumption that $x_n\in \Omega_j$, which implies that $x$ has at most $j$ nonzero components. This contradiction proves that $\Omega_j$ is closed.
\end{proof}

The next result translates the openness of $A_d$ to its preimage set $B_d$.

\begin{proposition}\label{Open-Bd}
If $M$ is a $d\times m$ full rank matrix, then the set $B_d$ is an open set in $\mathbb{R}^m$.
\end{proposition}
\begin{proof}
We prove this result by contradiction. Assume to the contrary that $B_d$ is not open. Then, there exists a point $\hat x\in B_d$ such that for all arbitrarily small $\epsilon>0$, the open balls
$$
B_o\left(\hat x, \frac{\epsilon}{\|M\|}\right):=\left\{x\in \mathbb{R}^m: \|x-\hat x\|_2<\frac{\epsilon}{\|M\|}\right\}
$$
of $\mathbb{R}^m$ are not completely contained in $B_d$. Here, $\|M\|$ denotes the spectral norm of the matrix $M$ induced by the Euclidean norm and $\|M\|>0$ is guaranteed by the hypothesis that $M$ is of full rank. Hence, for each $\epsilon$, there exists $x_\epsilon\in B_o(\hat x, \epsilon/\|M\|)$ such that $x_\epsilon\notin B_d$. Let $y_\epsilon:=Mx_\epsilon$. By Proposition \ref{A-B}, we have that $y_\epsilon\notin A_d$. Moreover, we let $\hat y:=M\hat x$. Since $\hat x\in B_d$, we clearly have that $\hat y\in A_d$. Therefore, we obtain that
$$
\|y_\epsilon-\hat y\|_2=\|M(x_\epsilon-\hat x)\|_2\leq \|M\|\|x_\epsilon-\hat x\|_2<\epsilon.
$$
This implies that the set $A_d$ is not open and contradicts part (i) of Proposition \ref{Closedness}. Therefore, $B_d$ is an open set in $\mathbb{R}^m$.
\end{proof}

It is clear that the $\ell_0$ norm is not a continuous function in the sense that
the condition $\|x_n-\hat x\|_2\to 0$ does not guarantee that $\|x_n\|_0\to \|\hat x\|_0$. To see this, we consider the sequence $x_n:=(1/n, 1/n, \dots, 1/n)$ and $\hat x=0$, the zero vector. Clearly, we have that $\|\hat x\|_0=0$ and
$$
\|x_n-\hat{x}\|_2=\sqrt{d}/n\to 0, \ \ \mbox{as}\ \ n\to\infty,
$$
but for all $n\in \mathbb{N}$, we find that $\|x_n\|_0=d$, which does not tend to zero.

It is important to understand how the sparsity of a vector in $\mathbb{R}^d$ influences the sparsity of vectors in its neighbourhood.  
To this end, for a given index set $\mathcal{I}\subseteq\mathbb{N}_d$
we define a subspace of $\mathbb{R}^d$ by letting
\begin{eqnarray}
\mathcal{C}_{\mathcal{I}}~:=\left\{x\in\mathbb{R}^d:\
S(x)\subseteq\mathcal{I}\right\},
\label{def:C_I}
\end{eqnarray}
where $S(x)$ denotes the support of $x\in \mathbb{R}^d$, that is, $$
S(x):=\{i\in \mathbb{N}_d: x_i\neq 0\},\ \ \mbox{for}\ \ x\in \mathbb{R}^d.
$$
Clearly, $\mathcal{C}_{\mathcal{I}}$ is convex. It is convenient to define the set
\begin{eqnarray}
\partial\mathcal{C}_{\mathcal{I}}:=\left\{x\in\mathbb{R}^d:\
S(x)=\mathcal{I}\right\}.
\label{def:bar-C_I}
\end{eqnarray}
We first establish a technical lemma.

\begin{lemma}\label{Distance}
If for some $\ell\in \mathbb{Z}_{d+1}$, $\hat x\in A_\ell$, then
\begin{equation}\label{distance-from-hatx}
    {\rm dist} \left(\hat x, A_\ell\setminus \partial\mathcal{C}_{S(\hat x)}\right)>0,
\end{equation}
where
$$
{\rm dist} (x, A):=\min\{\|x-z\|_2: z\in A\}.
$$
\end{lemma}
\begin{proof}
Since $\hat x\in A_\ell$, we may assume that $S(\hat x)=\{k_1, k_2, \dots, k_\ell\}$, where $1\leq k_1< k_2<\cdots<k_\ell\leq d$. It follows that $\hat x$ may be represented as
$$
\hat x=\sum_{j\in\mathbb{N}_\ell}\hat x_{k_j} e_{k_j}, \ \ \mbox{for some}\ \ \hat x_{k_j}\in\mathbb{R}\setminus\{0\}, \ j\in\mathbb{N}_\ell.
$$
For all $z\in  A_\ell\setminus \partial\mathcal{C}_{S(\hat x)}$, we have that $S(\hat x)\neq S(z)$, and there exist integers $k'_j$, $j\in\mathbb{N}_\ell$, with $1\leq k'_1<k'_2<\cdots<k'_\ell\leq d$ such that
$$
z=\sum_{j\in\mathbb{N}_\ell}z_{k'_j} e_{k'_j}, \ \ \mbox{for some}\ \ z_{k'_j}\in\mathbb{R}\setminus\{0\}, \ j\in\mathbb{N}_\ell.
$$
Hence, there exists some $k_j\in S(\hat x)$ but $k_j\notin S(z)$. This fact together with the above representations of $\hat x$ and $z$ implies that for all $z\in  A_\ell\setminus \partial\mathcal{C}_{S(\hat x)}$,
$$
\|x-z\|_2\geq |\hat x_{k_j}|>0.
$$
From this we conclude that
$$
{\rm dist} \left(\hat x, A_\ell\setminus \partial\mathcal{C}_{S(\hat x)}\right)\geq |\hat x_{k_j}|>0,
$$
which completes the proof of this lemma.
\end{proof}

With the help of Lemma \ref{Distance}, we prove the following proposition.

\begin{proposition}\label{zeroNormCompare}
The following statements hold true:

(i) If for some $\ell\in \mathbb{Z}_{d+1}$, $\hat x\in A_\ell$, then there exists $\delta>0$ such that for all $x\in B(\hat x, \delta)$, there holds $x\in \cup_{j=\ell}^d A_j$, that is, $\|x\|_0\geq \|\hat x\|_0$.

(ii) If for some $\ell\in \mathbb{Z}_{d+1}$, $\hat x\in A_\ell$, then there exists $\delta>0$ such that for all $x \in B(\hat x, \delta)\setminus\mathcal{C}_{\mathcal{I}}$ with $\mathcal{I}:=S(\hat x)$, there holds $x\in A_j$, for some $j\geq \ell+1$, that is, $\|x\|_0\geq \|\hat x\|_0+1$.
\end{proposition}
\begin{proof}
(i) We prove this assertion by contradiction. Assume to the contrary that the statement is not true. Then, for any $\delta>0$, there exists $x_\delta\in B(\hat x, \delta)$ such that $x_\delta\in \Omega_{\ell-1}$. By Item (ii) of Proposition \ref{Closedness}, the set $\Omega_{\ell-1}$
is closed. This implies that $\hat x\in \Omega_{\ell-1}$, which contradicts the assumption that $\hat x\in A_{\ell}$. Therefore, there exists $\delta>0$ such that for all $x\in B(\hat x, \delta)$, $x\in A_j$, for some $j\geq \ell$. This further implies that $\|x\|_0\geq \|\hat x\|_0$.

(ii) By Lemma \ref{Distance}, we may choose $\delta_0>0$ such that 
$$
{\rm dist} \left(\hat x, A_\ell\setminus \partial\mathcal{C}_{{S}(\hat x)}\right)>\delta_0.
$$
By Item (i) of this proposition, there exists a $\delta$ with $\delta_0>\delta>0$ such that for all $x\in B(\hat x, \delta)$, we have that $x\in A_j$, for some $j\geq \ell$. Hence, for this positive number $\delta$, there holds
$$
(B(\hat x, \delta)\setminus \partial\mathcal{C}_{{S}(\hat x)})\cap A_\ell=\emptyset.
$$
It follows that for all $x \in B(\hat x, \delta)\setminus\mathcal{C}_{{S}(\hat x)}$, there holds $x\in A_j$, for some $j\geq \ell+1$. Consequently, for all $x \in B(\hat x, \delta)\setminus\mathcal{C}_{{S}(\hat x)}$, we have that $\|x\|_0\geq \|\hat x\|_0+1$.
\end{proof}

\section{Sparsity Regularization under a Transform}

In this section we consider the minimization problem of sparsity regularization under a transform. The rationale for considering the regularization problem \eqref{regularization} with the function $f$ having the form \eqref{main-function-transform} is that $g$ has a global minimizer but it may not be sparse under the transform. We then impose the regularization term. By choosing the parameter $\lambda$ appropriately, we seek a global minimizer of $f$ having sparsity of a prescribed level and close to the global minimizer of $g$.  Specifically, we intend to understand how choices of the regularization parameter $\lambda$ lead to sparsity (under the transform $M$) of a global minimizer of the function $f$ defined by \eqref{main-function-transform} when $M$ is a real $d\times m$ matrix of full rank.

We first comment on a connection between the sets $A_\ell$ defined by equation \eqref{setA} and the $\ell_0$ norm. By the definition of the $\ell_0$ norm, for any $x\in \mathbb{R}^d$ we have that
\begin{equation}\label{l0-norm}
\|x\|_0=\ell, \ \ \mbox{if} \ \ x\in A_\ell, \ \ \mbox{for some}\ \ \ell\in
\mathbb{Z}_{d+1}.
\end{equation}
Formula \eqref{l0-norm} can simplify the function $f$ defined by \eqref{main-function-transform} on each set $B_\ell$ and provides a key to understand the solution of the related regularization problem \eqref{regularization}. In fact, by employing formula \eqref{l0-norm} and the partition $B_j$, $j\in \mathbb{Z}_{d+1}$, of $\mathbb{R}^m$, connected with the sets $A_\ell$, $\ell\in
\mathbb{Z}_{d+1}$ via Proposition \ref{A-B}, we have an alternative representation of function $f$ defined by \eqref{main-function-transform}. Namely,
\begin{equation}\label{Alternative2}
f(x)=g(x)+\lambda\ell, \ \ \mbox{for all}\ \ x\in B_\ell, \ \ell \in \mathbb{Z}_{d+1}.
\end{equation}

Geometric interpretation of the function $f$ defined by \eqref{main-function-transform} provides insights to sparsity of a global minimizer of $f$ under the transform $M$.
By adding the regularization term $\lambda\|M\cdot\|_0$ to $g$ results in lifting the graph of $g$ according to the sparsity in the range of $M$. In other words, the regularization term $\lambda\|M\cdot\|_0$ terraces the graph of function $g$. Specifically, the values of function $g$ that stay unchanged are $g(x)$ for all $x\in B_0$ (that is, in the null space of $M$) and every other values $g(x)$ are lifted according to which set $B_j$ the points $x$ belong to. For example, for all $x\in B_1$, the values $g(x)$ are lifted to $g(x)+\lambda$. In general, for all $x\in B_j$, the values $g(x)$ are lifted to $g(x)+\lambda j$, for $j\in \mathbb{Z}_{d+1}$. On the highest level of the terraces are $g(x)+\lambda d$, for all $x\in B_d$, where all components of $Mx$ are nonzero. Hence, by changing the value of the parameter $\lambda$, the landscape of the graph of the associated function $f$ is changed and accordingly the sparsity of the global minimizer of $f$ is changed. For instance, if the most sparse global minimizer is desired (that is, a point in the null space of $M$ as a global minimizer of $f$), then a value of $\lambda$ is chosen so that the function values $g(x)+\lambda j$, for all $x\in B_j$, $j\in \mathbb{N}_d$, are greater than the value $g(\hat x)$, where $\hat x$ is in the null space of the matrix $M$. This understanding is a key to guide for choices of the parameter $\lambda$.


When a global minimizer of $f$ that are most sparse is desired, we have the parameter choice strategy described in the next theorem.

\begin{theorem}\label{Highest}
Let $x^*\in \mathbb{R}^m$ be a global minimizer of $g$ and $x_0\in B_0$ be a minimizer of $g$ on $B_0$. If the parameter $\lambda$ is chosen to satisfy
\begin{equation}\label{x-zero}
\lambda\geq g(x_0)-g(x^*),
\end{equation}
then $x_0$ is a global minimizer of $f$ in the space $\mathbb{R}^m$, $Mx_0$ has sparsity of level $0$ and the minimum value of $f$ is given by $g(x_0)$.
\end{theorem}
\begin{proof}
Let $x\in \mathbb{R}^m$ be an arbitrarily fixed vector. We make use of the partition $B_j$, $j\in\mathbb{Z}_{d+1}$, of $\mathbb{R}^m$ to conclude that there exists $j'\in \mathbb{Z}_{d+1}$ such that $x\in B_{j'}$.
We consider two cases: $j'=0$ and $j'\in \mathbb{N}_d$.

In the case when $j'=0$, we have that $x\in B_0$ and $Mx=0$. In this case, by employing equation \eqref{Alternative2}, from the assumption that $x_0\in B_0$ is a minimizer of $g$ on $B_0$ we get that
$$
f(x)=g(x)\geq g(x_0)=f(x_0).
$$

Next, we consider the case when $j'\in\mathbb{N}_d$. In this case, we have that $x\in B_{j'}$, that is, $Mx\in A_{j'}$. Once again,
we employ equation \eqref{Alternative2} to obtain that
$$
f(x)=g(x)+\lambda j'\geq g(x)+\lambda.
$$
Combining this inequality with the assumption on $x^*$ and condition \eqref{x-zero}, we obtain that
$$
f(x)\geq g(x)+\lambda\geq g(x^*)+\lambda \geq g(x_0)=f(x_0).
$$

In both of the cases, we have shown that
$
f(x)\geq f(x_0)\ \  \mbox{for all}\ \ x\in \mathbb{R}^m.
$
Therefore, $x_0$ is a global minimizer of $f$ on the space $\mathbb{R}^m$.
\end{proof}

Theorem \ref{Highest} also provides the error bound of the regularized global minimum value $f(x_0)$ from the original global minimum value $g(x^*)$. Namely,
$$
0\leq f(x_0)-g(x^*)\leq\lambda.
$$
In general, the error is not equal to zero unless the global minimizer $x^*$ of $g$ is in $B_0$, in which case sparse  regularization is not necessary.

\begin{figure}[htbp]
 \centering
 \begin{tabular}{cc}
 \includegraphics[width=6.0cm,height=5.5cm]{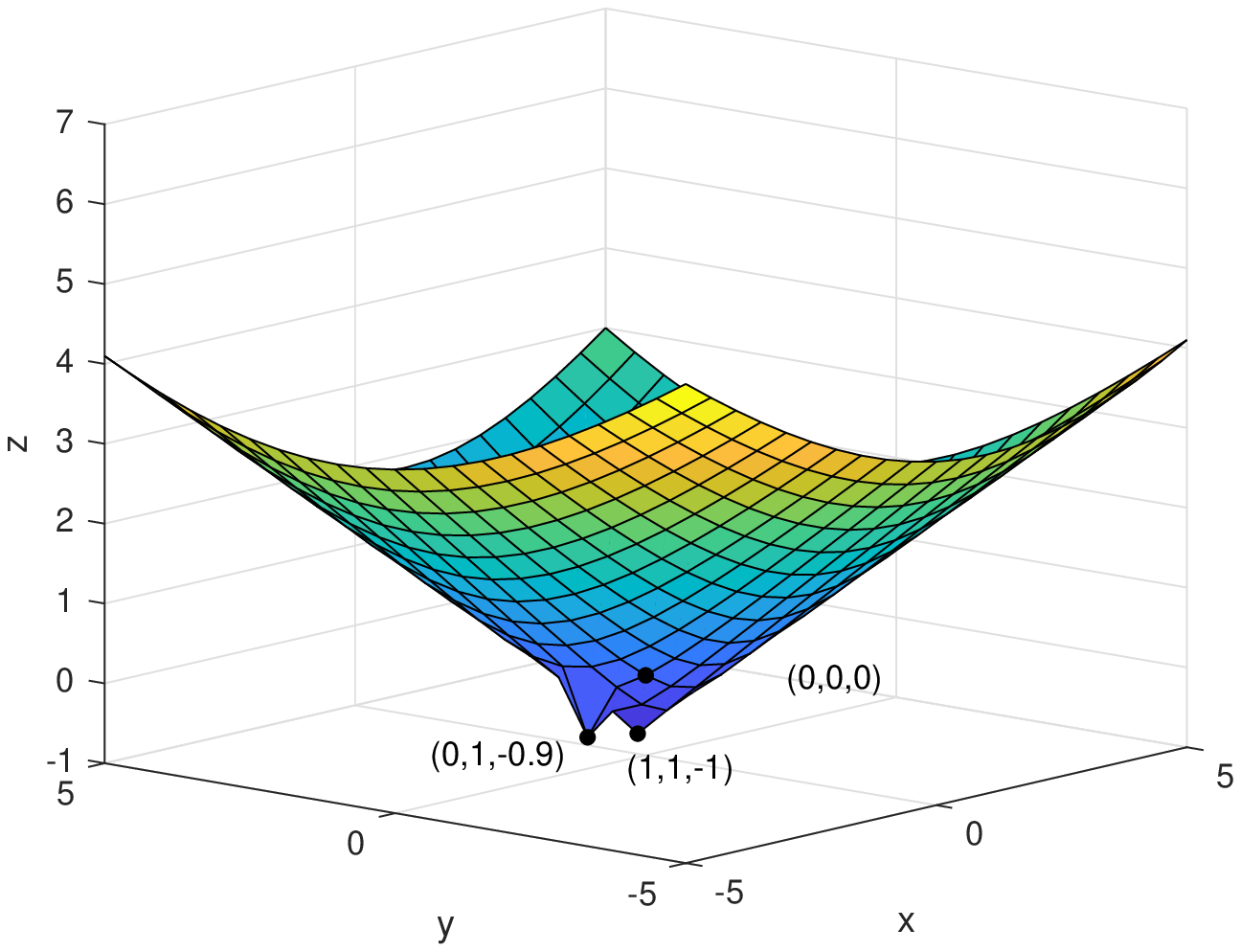} 
 & \includegraphics[width=6.0cm,height=5.5cm]{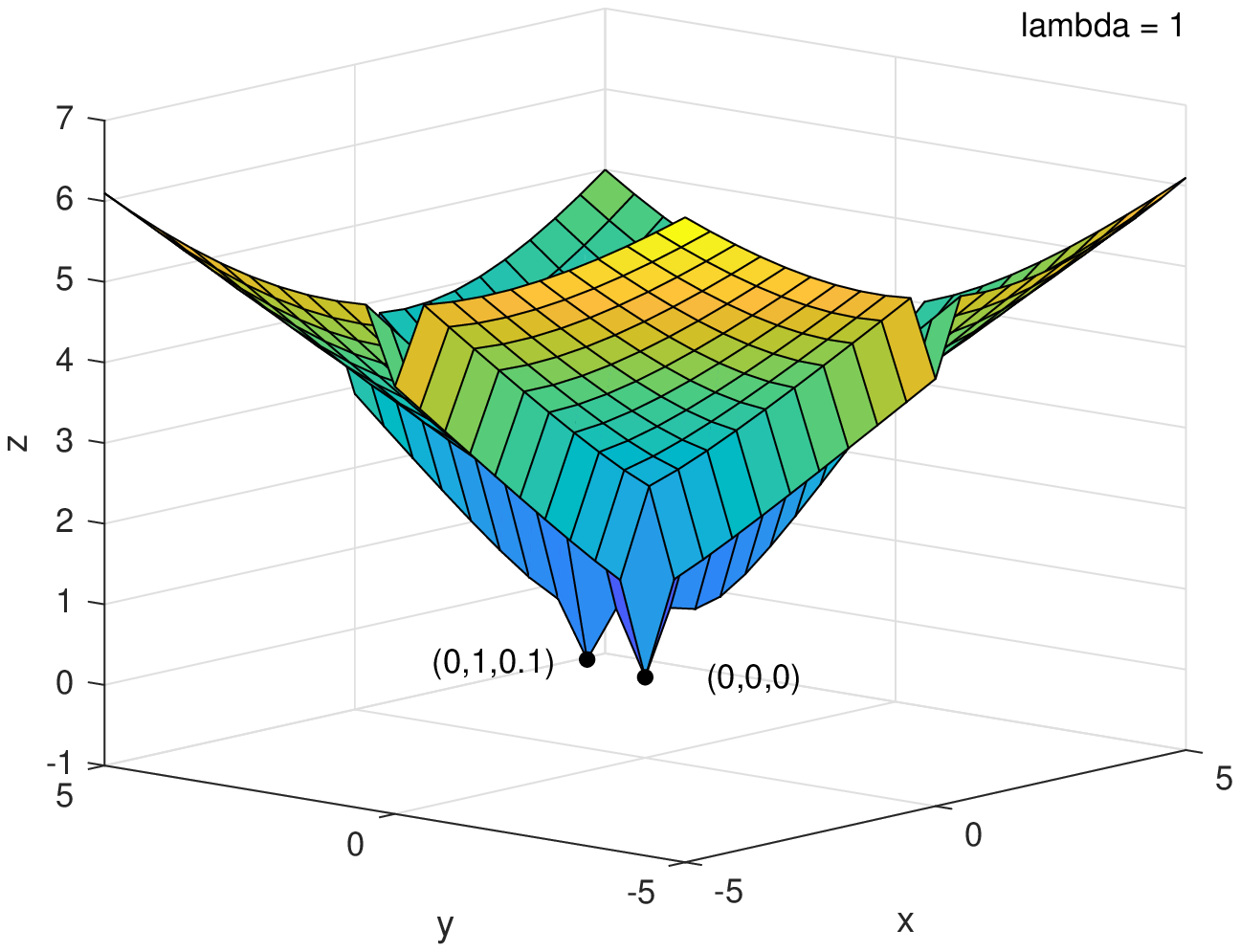}\\
 \scriptsize{(a)}&\scriptsize{(b)}
 \end{tabular}
 \caption{ (a) Graph of function $g$; (b) Graph of the regularized function $f$ with $\lambda=1$.}
 \label{fig:geometric}
\end{figure}

We illustrate Theorem \ref{Highest} by a simple example in $\mathbb{R}^2$. To this end, we consider a non-convex function defined by
\begin{equation}\label{exampleg}
g({x}):=
\left\{
  \begin{array}{ll}
    \frac{\sqrt{2}}{2}\|{x}-\mathbf{1}\|_2-1,\ \ & {x}\neq(0,1), \\
    -0.9, & {x}=(0,1),
  \end{array}
\right.
\end{equation}
where ${\bf 1}:=(1,1)$.
Clearly, as shown in Fig. \ref{fig:geometric} (a)  $g$ has a unique, {\it non-sparse} global minimizer $x^*=(1,1)$ and the minimum value $g(x^*)=-1$. In this example, we choose $M=I$. The value $g(x_0)=0$, where $x_0:=(0,0)$. According to Theorem \ref{Highest}, we choose 
$$
\lambda=g(x_0)-g(x^*)=1.
$$
The regularized non-convex function $f$ defined by \eqref{main-function} with $g$ in the form \eqref{exampleg} and $\lambda=1$ is shown in Fig. \ref{fig:geometric} (b). Note that ${x}_0:=(0,0)$ is the {\it sparse} global minimizer of $f$ and the minimum value $f({x}_0)=0$. Also, the error between the regularized global minimizer and the original global minimizer is given by $f(x_0)-g(x^*)=1=\lambda$. These figures 
illustrate how the regularization term $\lambda\|\cdot\|_0$ terraces 
the graph of $g$.


In the case that we wish to reduce the error of the regularized global minimum value, we may choose not to demand the most sparsity (under the transform $M$). For this reason, we consider a choice of the parameter $\lambda$ with which the function $f$ defined by equation \eqref{main-function-transform} has a global minimizer having sparsity (under the transform $M$) of a prescribed level. We present a parameter choice strategy in the next theorem.
For this purpose, we find it convenient to define another sequence of sets. For all $\ell\in \mathbb{Z}_{d+1}$, we let
\begin{equation}\label{Gamma}
\Gamma_\ell:=\bigcup_{j\in \mathbb{Z}_{\ell+1}} B_j.
\end{equation}
It can be readily verified that for each $\ell\in \mathbb{Z}_{d+1}$, the set $\Gamma_\ell$ is the preimage set of $\Omega_\ell$ under the transform $M$. Moreover, by \eqref{partitionB}, we have that
\begin{equation}\label{Gamma-d}
\Gamma_d=\mathbb{R}^m.
\end{equation}

\begin{theorem}\label{choice-of-lambda2}
Let $x^*\in \mathbb{R}^m$ be a global minimizer of $g$, $x'\in \Gamma_\ell$ be a minimizer of $g$ on $\Gamma_\ell$ for some $\ell\in \mathbb{N}_d$, and $x_{j}\in \Gamma_{j}$ be a minimizer of $g$ on $\Gamma_{j}$, for all $j\in \mathbb{Z}_\ell$.
Suppose that
\begin{equation}\label{conditions-on-min2}
g(x')-g(x^*)\leq \frac{1}{\ell-j}[g(x_j)-g(x')], \ \ \mbox{for all}\ \ j\in \mathbb{Z}_\ell.
\end{equation}
If the parameter $\lambda$ is chosen to satisfy the conditions
\begin{equation}\label{conditions-on-lambda2}
g(x')-g(x^*)\leq \lambda\leq\frac{1}{\ell-j}[g(x_{j})-g(x')],\ \ \mbox{for all}\ \ j\in \mathbb{Z}_\ell,
\end{equation}
then $x'$ is a global minimizer of $f$ on $\mathbb{R}^m$, $Mx'$ has sparsity of level $\ell'$ with $\ell'\leq \ell$ and the global minimum value of $f$ on $\mathbb{R}^m$ is given by $f(x')=g(x')+\lambda \ell'$.
\end{theorem}
\begin{proof}
We shall verify that
\begin{equation}\label{G-min2}
f(x')\leq f(x), \ \ \mbox{for all}\ \ x\in \mathbb{R}^m.
\end{equation}
Since $x'\in \Gamma_\ell$ is a minimizer of $g$ on $\Gamma_\ell$, we use the definition \eqref{Gamma} of $\Gamma_\ell$ to consider cases when $x'\in B_j$ for all $j\in\mathbb{Z}_{\ell+1}$.

We consider the first case when $x'\in B_0$, that is, $Mx'\in A_0$. In this case, we shall show \eqref{G-min2} with $Mx'=0$ by verifying it for all $x\in B_j$, for all $j\in \mathbb{Z}_{d+1}$. To this end, we consider two subcases according to the index $j\in\mathbb{Z}_{d+1}$: (1) $j\leq \ell$ and (2) $j>\ell$. In subcase (1), we consider $x\in B_j$, $j\in \mathbb{Z}_{\ell+1}$. By the definition of the sets $B_j$, we have that $Mx\in A_j$, which implies
$$
\|Mx\|_0=j,\ \ \mbox{if}\ \ x\in B_j,\ \ \mbox{for all}\ \ j\in \mathbb{Z}_{\ell+1}.
$$
Since $x'$ is a minimizer of $g$ on $\Gamma_\ell$ and $\|Mx'\|_0=0$, we have for all $x\in B_j$, $j\in \mathbb{Z}_{\ell+1}$, that
$$
f(x')=g(x')\leq g(x)\leq g(x)+\lambda j=f(x).
$$
In subcase (2), we consider $x\in B_j$, $j=\ell+1, \ell+2, \dots, d$, where we also have that $Mx\in A_j$ and thus,
$\|Mx\|_0=j$.
By employing the first inequality of \eqref{conditions-on-lambda2} and the fact that $x^*$ is a global minimizer of $g$, for all $x\in B_j$, $j=\ell+1, \ell+2, \dots, d$, we obtain that
$$
f(x')=g(x')\leq g(x^*)+\lambda \leq g(x)+\lambda\leq g(x)+\lambda j=f(x).
$$
This proves the first case of  \eqref{G-min2}.

We next consider the second case when $x'\in B_k$, for $k\in \mathbb{N}_\ell$ and show \eqref{G-min2} with $\|Mx'\|_0=k$ by verifying \eqref{G-min2} for all $x\in B_j$, for $j\in \mathbb{Z}_{d+1}$. To this end, we consider three subcases according to the index $j\in\mathbb{Z}_{d+1}$: (1) $j<\ell$, (2) $j=\ell$, (3) $j>\ell$. In subcase (1) for which $j<\ell$, we consider $x\in B_j$, $j\in \mathbb{Z}_{\ell}$. Clearly, we have that $\|Mx\|_0=j$. Thus, by using the second inequality of \eqref{conditions-on-lambda2} and the hypothesis that $x_j\in \Gamma_j$ is a minimizer of $g$ on $\Gamma_j$, we observe that
$$
g(x')+\lambda \ell\leq g(x_j)+ \lambda j\leq g(x)+ \lambda j=f(x).
$$
This together with the fact $k\leq \ell$ implies that
$$
f(x')=g(x')+\lambda k\leq g(x')+\lambda \ell\leq f(x), \ \ \mbox{for all}\ \ x\in B_j, \ \ \mbox{for all}\ \ j\in \mathbb{Z}_{\ell}.
$$
In subcase (2) for which $j=\ell$, we consider $x\in B_\ell$. Since $x'$ is a minimizer of $g$ on $\Gamma_\ell$, which contains $B_\ell$ as a subset, and $k\leq \ell$, we find that
$$
f(x')=g(x')+\lambda k\leq g(x)+\lambda k\leq g(x)+\lambda \ell=f(x), \ \ \mbox{for all}\ \ x\in B_\ell.
$$
In subcase (3) for which $j>\ell$, we consider $x\in B_j$, for $j=\ell+1, \ell+2, \dots, d$. By using the first inequality of \eqref{conditions-on-lambda2} and again the fact that $x^*$ is a global minimizer of $g$, we derive for all $x\in B_j$, for all $j=\ell+1, \ell+2, \dots, d$ that
$$
f(x')=g(x')+\lambda k\leq g(x^*)+\lambda(k+1)\leq g(x)+\lambda j=f(x).
$$
That is, \eqref{G-min2} holds true for the second case.

Summarizing the above verification, we conclude that \eqref{G-min2} holds true for all cases and thus, $x'$ is a global minimizer of the function $f$ on the space $\mathbb{R}^m$.
\end{proof}

A comment on the parameter choice \eqref{conditions-on-lambda2} in Theorem \ref{choice-of-lambda2} is in order. To be able to choose such a parameter, it requires that condition \eqref{conditions-on-min2} is satisfied.
This hypothesis is indeed needed to ensure that the choice \eqref{conditions-on-lambda2} of the parameter $\lambda$ is feasible. The hypothesis \eqref{conditions-on-min2} is equivalent to the following conditions that
\begin{equation}\label{average}
g(x')\leq \frac{g(x_j)+(\ell-j)g(x^*)}{\ell-j+1}, \ \ \mbox{for all} \ \ j\in \mathbb{Z}_\ell.
\end{equation}
The right-hand-side of the inequality \eqref{average} is a weighted average of $g(x_j)$ and $g(x^*)$.
By the definition of $x'$, $x_j$ and $x^*$, it is clear that
$$
g(x^*)\leq g(x')\leq g(x_j), \ \ \mbox{for all} \ \ j\in \mathbb{Z}_\ell.
$$
This shows that the hypothesis \eqref{conditions-on-min2} of Theorem \ref{choice-of-lambda2} is reasonable. Condition \eqref{conditions-on-lambda2} also reveals that the error of the regularized global minimum value $f(x')$ approximating the original global minimum value $g(x^*)$ is bounded by the value of the regularization parameter so chosen. That is,
$$
0\leq f(x')-g(x^*)\leq \lambda.
$$

We next illustrate the result of  Theorem \ref{choice-of-lambda2} by presenting a corollary of Theorem \ref{choice-of-lambda2} for the special case when $\ell=1$. The corollary gives a choice of the parameter $\lambda$ which guarantees that a global minimizer of $f$ has sparsity of level $1$ under the transform $M$. That is, the corresponding function $f$ has a global minimizer having at most one nonzero component under the transform $M$.

\begin{coro}
Let $x^*\in \mathbb{R}^d$ be a global minimizer of $g$ on the space $\mathbb{R}^m$, $x'\in \Gamma_1$ be a minimizer of $g$ on $\Gamma_1$ and $x_0\in B_0$ be a minimizer of $g$ on $B_0$.
If the parameter $\lambda$ is chosen to satisfy the condition
\begin{equation}\label{conditions-on-lambda-ell=12}
g(x')-g(x^*)\leq \lambda\leq g(x_0)-g(x'),
\end{equation}
then $x'$ is a global minimizer of $f$ on the space $\mathbb{R}^m$.
\end{coro}
\begin{proof}
This result is obtained by specializing Theorem \ref{choice-of-lambda2} to the special case when $\ell=1$.
\end{proof}

Clearly, as we have discussed earlier, for the choice \eqref{conditions-on-lambda-ell=12} of the parameter $\lambda$ to be feasible, we need to require that
\begin{equation}\label{conditions-on-min-ell=12}
g(x')\leq \frac{1}{2}[g(x^*)+g(x_0)].
\end{equation}
That is, $g(x')$ is less than or equal to the average of $g(x^*)$ and $g(x_0)$.

The next theorem connects a global minimizer of $g$ with a global (or local) minimizer of $f$. To this end, we recall the definition of a local minimizer of a non-convex function. We first define a closed ball centered at $x^*\in\mathbb{R}^d$ with radius $\delta>0$ by
$$
B(x^*, \delta):=\{x\in \mathbb{R}^d: \|x-x^*\|_2\leq \delta\}.
$$
A vector $x^*\in \mathbb{R}^d$ is called a local minimizer of $f$, if there exists a $\delta>0$ such that
$$
f(x^*)
\leq f(x),
\ \ \mbox{for all }\ \ x\in B(x^*, \delta).
$$

\begin{theorem}
Let $x^*\in \mathbb{R}^m$ be a global minimizer of $g$.

(i) If $Mx^*=0$, then $x^*$ is a global minimizer of $f$ on  $\mathbb{R}^m$.

(ii) If for some $\ell\in \mathbb{N}_d$, $x^*\in B_\ell$ and if for a minimizer $x_j\in B_j$ of $g$ on $B_j$ for all $j\in \mathbb{Z}_\ell$, the parameter $\lambda$ is chosen to satisfy
\begin{equation}\label{conditionii}
0\leq \lambda \leq \frac{1}{\ell-j}[g(x_j)- g(x^*)],\ \ \mbox{for all}\ \ j\in
\mathbb{Z}_\ell,
\end{equation}
then $x^*$ is a global minimizer of $f$ on  $\mathbb{R}^m$.

(iii) If $x^*\in B_d$,
then $x^*$ is a local minimizer of $f$.

(iv) If $x^*\in B_d$, and for some $j\in \mathbb{Z}_{d+1}$ and for some $\tilde x\in B_j$,
\begin{equation}\label{conditioniv}
g(x^*)+\lambda(d-j)>g(\tilde x),
\end{equation}
then $x^*$ is a local minimizer of $f$ but not a global minimizer of $f$ on $\mathbb{R}^m$, and global minimizers (if exist) of $f$ on $\mathbb{R}^m$ have sparsity of level at least $d-1$.
\end{theorem}
\begin{proof}
Since $x^*\in \mathbb{R}^d$ is a global minimizer of $g$, we have
that
\begin{equation}\label{g-min-d2}
g(x^*)\leq g(x), \ \ \mbox{for all} \ \ x\in \mathbb{R}^m.
\end{equation}
For both Items (i) and (ii), we shall show that
\begin{equation}\label{f-min-d2}
f(x^*)\leq f(x), \ \ \mbox{for all} \ \ x\in \mathbb{R}^m.
\end{equation}

(i) If $Mx^*=0$, by the definition of $\|\cdot\|_0$, we have that $\|Mx^*\|_0=0$ and thus,
$$
f(x^*)=g(x^*)+\|Mx^*\|_0=g(x^*).
$$
Consequently, according to condition
\eqref{g-min-d2}, we obtain that
$$
f(x^*)=g(x^*)\leq g(x)\leq g(x)+\lambda\ell=f(x), \ \ \mbox{for} \ \ x\in
B_\ell, \ \ \mbox{for all}\ \ \ell\in \mathbb{Z}_{d+1}.
$$
This confirms that \eqref{f-min-d2} holds true.

(ii) Since for some $\ell\in \mathbb{N}_d$, $x^*\in B_\ell$ and satisfies condition \eqref{conditionii}, and since $x_j\in B_j$ is a minimizer of $g$ on $B_j$ for all $j\in \mathbb{Z}_\ell$, we have that
$$
f(x^*)=g(x^*)+\lambda\ell\leq g(x_j)+\lambda j\leq g(x)+\lambda j=f(x), \ \ \mbox{for all}
\ \ x\in B_j, \ \ j\in
\mathbb{Z}_\ell.
$$
Moreover, we have that
$$
f(x^*)=g(x^*)+\lambda\ell\leq g(x)+\lambda j=f(x), \ \ \mbox{for all} \ \ x\in
B_{j}, \ \ j=\ell, \ell+1, \dots, d.
$$
Hence, \eqref{f-min-d2} is satisfied.

(iii) Since $x^*\in B_d$ is a global minimizer of $g$ on $\mathbb{R}^m$, we have for all $x\in B_d$ that
$$
f(x^*)=g(x^*)+\lambda d\leq g(x)+\lambda d=f(x).
$$
That is,
\begin{equation}\label{local2}
f(x^*)\leq f(x), \ \ \mbox{for all}\ \ x\in B_d.
\end{equation}
According to Proposition \ref{Open-Bd}, $B_d$ is an open set. Thus, $x^*$ is an interior point of $B_d$. This ensures that there exists a $\delta>0$ such that $B(x^*, \delta)$ is contained in $B_d$.  Therefore, from inequality \eqref{local2} we conclude that $f(x^*)\leq f(x)$ for all $x\in B(x^*, \delta)$, and thus, $x^*$ is a local minimizer of $f$ on $\mathbb{R}^m$.

(iv) By (iii), we have known that in this case, $x^*\in B_d$ is a local minimizer of $f$. We next show that $x^*$ is not a global minimizer of $f$ on $\mathbb{R}^m$. By \eqref{conditioniv}, we have that for some $j\in \mathbb{Z}_{d+1}$ and for some $\tilde x\in B_j$,
$$
g(x^*)+\lambda d>g(\tilde x)+\lambda j.
$$
This together with the fact $x^*\in B_d$ ensures that
$$
f(x^*)=g(x^*)+\lambda d >g(\tilde x)+\lambda j= f(\tilde x).
$$
This implies that $x^*$ is not a global minimizer of $f$.

Finally, we prove that there is a global minimizer of $f$ on $\mathbb{R}^m$ having sparsity of level at least $d-1$. From \eqref{local2}, we know that $x^*$ is a minimizer of $f$ on $B_d$. Note that by \eqref{Gamma-d}, there holds $\mathbb{R}^m\setminus B_d=\Gamma_{d-1}$. Hence, the fact established earlier that $x^*$ is not a global minimizer of $f$ implies that global minimizers (if exist) of $f$ must occur at a point in $\Gamma_{d-1}$. By the definition of $\Gamma_{d-1}$ such a global minimizer has sparsity of level at least $d-1$.
\end{proof}

In the next theorem, we prove necessary conditions of a global minimizer of $f$.

\begin{theorem}\label{necessary-condition2}
Let $x^*\in\mathbb{R}^m$ be a global minimizer of $f$ on $\mathbb{R}^m$.

(i) If $x^*\in B_\ell$ for some $\ell\in \mathbb{Z}_{d+1}$, then $x^*$ is a minimizer of $g$ on $\Gamma_\ell$.

(ii) If $x^*\in\mathbb{R}^m$ is not a global minimizer of $g$ on $\mathbb{R}^m$, then $x^*\in \Gamma_{d-1}$.
\end{theorem}
\begin{proof}
(i) Since $x^*\in B_\ell$ for some $\ell\in \mathbb{Z}_{d+1}$ and it is
a global minimizer of $f$ on $\mathbb{R}^m$,  we have that
$$
g(x^*)+\lambda \ell=f(x^*)\leq f(x)=g(x)+\lambda j, \ \ \mbox{for all} \ \ x\in B_j, \ \ j\in \mathbb{Z}_{\ell+1}.
$$
It follows that
$$
g(x^*)+\lambda(\ell-j)\leq g(x), \ \ \mbox{for all} \ \ x\in B_j, \ \ j\in \mathbb{Z}_{\ell+1}.
$$
Using this inequality and noting that
$$
\lambda (\ell-j)\geq 0, \ \ \mbox{for all}\ \ j\in \mathbb{Z}_{\ell+1},
$$
we have that
$$
g(x^*)\leq g(x^*)+\lambda(\ell-j)\leq g(x), \ \ \mbox{for all} \ \ x\in B_j,
\ \ j\in \mathbb{Z}_{\ell+1}.
$$
The above inequality together with the definition \eqref{Gamma} of the set $\Gamma_\ell$ ensures that $x^*\in B_\ell$ is a minimizer of $g$ on $\Gamma_\ell$.

(ii) We prove this assertion by contradiction. Assume to the contrary that $x^*\notin \Gamma_{d-1}$. Since $x^*\in\mathbb{R}^m$,
$$
\mathbb{R}^m=\Gamma_{d-1}\cup B_d
$$
and
$$
\Gamma_{d-1}\cap B_d=\emptyset,
$$
we must have that $x^*\in B_d$. By Statement (i) of this theorem with $\ell=d$, we conclude that $x^*$ is a minimizer of $g$ on $\Gamma_d$. Noting that $\Gamma_d=\mathbb{R}^m$, we confirm that $x^*$ is a global minimizer of $g$ on $\mathbb{R}^m$. This contradicts the hypothesis that $x^*$ is not a global minimizer of $g$ on $\mathbb{R}^m$. This contradiction ensures that $x^*\in \Gamma_{d-1}$.
\end{proof}

In Theorem \ref{necessary-condition2} (ii), we provide sufficient conditions which guarantee that a global minimizer of $f$ is sparse under the transform $M$.

The next result follows immediately from Theorem \ref{necessary-condition2} (ii).

\begin{coro}
If $x^*\in\mathbb{R}^m$ is a global minimizer of $f$ on $\mathbb{R}^m$,  then, either
$x^*\in \Gamma_{d-1}$ or $x^*$ is a global
minimizer of $g$ on $\mathbb{R}^m$.
\end{coro}

\section{Sparse Regularization in the Spacial Domain}

This section is devoted to presentation of special results for regularization having sparsity in the spacial domain. In this case, $d=m$ and $M=I$. All results in the last section can be restricted to the function defined by \eqref{main-function}. We leave the deviation of these results to the interested reader.

We consider in this section the optimization problem
\begin{equation}\label{opt:f}
\min\left\{f(x,y): (x,y) \in \mathbb{R}^d\times\mathbb{R}^{d'}\right\},
\end{equation}
where $f$ is defined by \eqref{def:f}, and present special results for this case. In this model, we seek sparsity for the variable $x$ only.
A typical example of optimization problem \eqref{opt:f} is approximately sparse regularization. In such cases, the function $f$ may take the following form
\begin{equation}\label{def:f2}
f(x,y):=\phi(y)+\mu\|x-Dy\|_2^2+\lambda\|x\|_0, \ \ (x,y)\in \mathbb{R}^d\times \mathbb{R}^{d'},
\end{equation}
or
\begin{equation}\label{def:f3}
f(x,y):=\phi(y)+\mu\|x-Dy\|_1+\lambda\|x\|_0, \ \ (x,y)\in \mathbb{R}^d\times \mathbb{R}^{d'},
\end{equation}
where $\phi$ is a convex function and $D$ is $d\times d'$ matrix. Form \eqref{def:f2} relates to regularization by the envelope of the $\ell_0$ norm and form  \eqref{def:f3} relates to regularization by the capped $\ell_1$ norm \cite{Jiang}. For specific examples of $f$, see \cite{WuXu2021} for inverting incomplete Fourier transform, \cite{Zeng18, Zeng19} for image/signal processing, \cite{Zheng19} for medical image reconstruction and  machine learning \cite{Lopez11, Louizos18, Pan, Wang2016}.

Employing the partition $A_j$, $j\in \mathbb{Z}_{d+1}$, of $\mathbb{R}^d$ and
the definition of $\|\cdot\|_0$, we have an alternative representation of function $f$:
\begin{equation}\label{Alternative}
f(x,y)=g(x,y)+\lambda\ell, \ \ \mbox{for all}\ \ x\in A_\ell, \ \ell \in \mathbb{Z}_{d+1} \ \ \mbox{and for all} \ \ y\in \mathbb{R}^{d'}.
\end{equation}
Clearly, adding the function $\lambda\|x\|_0$ to $g(x,y)$ results in lifting the graph of $g(x,y)$ according to the sparsity partition of $\mathbb{R}^d$ with respect to the variable $x$. In other words, the $\ell_0$ norm terraces the graph of function $g(x,y)$. Specifically, the only value that stays unchanged is $g(0,y)$ and every other value of $g(x,y)$ is lifted according to which set $A_j$ the point $x$ belongs to. For example, for $x\in A_1$, $g(x,y)$ is lifted to $g(x,y)+\lambda$. In general, for $x\in A_j$, $g(x,y)$ is lifted to $g(x,y)+\lambda j$, for $j\in \mathbb{Z}_{d+1}$. On the highest level of the terraces is $g(x,y)+\lambda d$, for all $x\in A_d$, where all components of $x$ are nonzero. Hence, by changing the value of the parameter $\lambda$, the landscape of the graph of the corresponding function $f$ is changed. Accordingly the sparsity of the global minimizer of $f$ is changed.

We first consider a choice of the parameter $\lambda$ with which the function $f$ defined by \eqref{def:f} has the most sparse minimizer in variable $x$. To this end, we assume that the function $g(x,y)$ has a global minimizer $(x^*,y^*)\in \mathbb{R}^d\times\mathbb{R}^{d'}$.

\begin{theorem}\label{most-sparse-set}
Let $(x^*,y^*)\in \mathbb{R}^d\times\mathbb{R}^{d'}$ be a global minimizer of $g$. If the parameter $\lambda$ is chosen to satisfy
\begin{equation}\label{zero}
\lambda\geq g(0,y_0)-g(x^*,y^*),
\end{equation}
for some $y_0\in \mathbb{R}^{d'}$,
then the pair $(0,y_0)\in \mathbb{R}^d\times\mathbb{R}^{d'}$ is a global minimizer of $f$ on $\mathbb{R}^d\times\mathbb{R}^{d'}$ and $g(0,y_0)$ is the minimum value of $f$ on $\mathbb{R}^d\times\mathbb{R}^{d'}$. Moreover, if the inequality \eqref{zero} becomes strict, then the pair $(0,y_0)\in \mathbb{R}^d\times\mathbb{R}^{d'}$ is the unique global minimizer of $f$ on $\mathbb{R}^d\times\mathbb{R}^{d'}$.
\end{theorem}
\begin{proof}
We consider an arbitrary $x\in \mathbb{R}^d\setminus\{0\}$ and use the sparsity partition of $\mathbb{R}^d$. There exists $j\in \mathbb{N}_d$ such that $x\in A_j$.
We employ equation \eqref{Alternative} to get
$$
f(x, y)=g(x, y)+\lambda j\geq g(x,y)+\lambda, \ \ \mbox{for all}\ \ y\in \mathbb{R}^{d'}.
$$
Combining this inequality with condition \eqref{zero}, we obtain that
\begin{equation}\label{f(x,y)>g(0,y*)}
    f(x,y)\geq g(x,y)+\lambda\geq g(x^*,y^*)+\lambda \geq g(0,y_0), \  \mbox{for all} \ x\in \mathbb{R}^d\setminus\{0\}\  \mbox{and all}\  y\in \mathbb{R}^{d'}.
\end{equation}
Moreover, by the definition of $f$, we have that
\begin{equation}\label{f(0,y*)}
f(0,y_0)=g(0,y_0)+\lambda \|0\|_0=g(0,y_0).
\end{equation}
Combining inequality \eqref{f(x,y)>g(0,y*)} and equation \eqref{f(0,y*)} yields that
$$
f(x,y)\geq f(0,y_0)\ \  \mbox{for all}\ \ x\in \mathbb{R}^d\setminus\{0\} \ \ \mbox{and for all}\ \ y\in \mathbb{R}^{d'}.
$$
Therefore, the pair $(0, y_0)\in \mathbb{R}^d\times \mathbb{R}^{d'}$ is a global minimizer of $f$. The uniqueness of the global minimizer of $f$ is guaranteed if a strict inequality for $\lambda$ is imposed.
\end{proof}

We next consider a choice of the parameter $\lambda$ with which the function $f(x,y)$ defined by  \eqref{def:f}  has a global minimizer with sparsity of a general level for the variable $x$.

\begin{theorem}\label{choice-of-lambda}
Let $(x^*,y^*)\in \mathbb{R}^d\times\mathbb{R}^{d'}$ be a global minimizer of $g$, for some $\ell\in \mathbb{N}_d$, $(x',y')\in \Omega_\ell\times\mathbb{R}^{d'}$ be a minimizer of $g$ on $\Omega_\ell\times\mathbb{R}^{d'}$ and $(x_j,y_j)\in \Omega_j\times\mathbb{R}^{d'}$ be a minimizer of $g$ on $\Omega_j\times\mathbb{R}^{d'}$, for all $j\in \mathbb{Z}_\ell$. Suppose that
\begin{equation}\label{conditions-on-min}
g(x',y')-g(x^*,y^*)\leq \frac{1}{\ell-j}[g(x_j,y_j)-g(x',y')], \ \ \mbox{for all}\ \ j\in \mathbb{Z}_\ell.
\end{equation}
If the parameter $\lambda$ is chosen to satisfy the conditions
\begin{equation}\label{conditions-on-lambda}
g(x',y')-g(x^*,y^*)\leq \lambda\leq\frac{1}{\ell-j}[g(x_j,y_j)-g(x',y')], \ \ \mbox{for all}\ \ j\in \mathbb{Z}_\ell,
\end{equation}
then $(x',y')$ is a global minimizer of $f$ on $\mathbb{R}^d\times\mathbb{R}^{d'}$.
\end{theorem}
\begin{proof}
It suffices to verify that
\begin{equation}\label{G-min}
f(x',y')\leq f(x,y), \ \ \mbox{for all}\ \ (x,y)\in \mathbb{R}^d\times\mathbb{R}^{d'}.
\end{equation}
Since $(x',y')\in \Omega_\ell\times\mathbb{R}^{d'}$ is a minimizer of $g$ on $\Omega_\ell\times\mathbb{R}^{d'}$, we use the definition of $\Omega_\ell$ to consider cases when $x'\in A_j$ for $j\in\mathbb{Z}_{\ell+1}$.

Step one: we consider the case when $x'\in A_0$, that is, $x'=0$. We now show \eqref{G-min} with $x'=0$ by verifying it for all $x\in A_j$, for $j\in \mathbb{Z}_{d+1}$, and for all $y\in\mathbb{R}^{d'}$. The case for $j=0$ is trivial and we consider other cases. Since $(x',y')\in \Omega_\ell\times\mathbb{R}^{d'}$ is a minimizer of $g$ on $\Omega_\ell\times\mathbb{R}^{d'}$ and $\|x'\|_0=0$, we have for all $x\in A_j$, $j=1,2,\dots, \ell$, and for all $y\in\mathbb{R}^{d'}$ that
$$
f(x',y')=g(x',y')\leq g(x,y)\leq g(x,y)+\lambda j=f(x,y).
$$
For all $x\in A_j$, $j=\ell+1, \ell+2, \dots, d$, and for all $y\in \mathbb{R}^{d'}$, by employing the first inequality of \eqref{conditions-on-lambda} and the assumption that $(x^*,y^*)$ is a global minimizer of $g$ on $\mathbb{R}^d\times\mathbb{R}^{d'}$, we obtain that
$$
f(x',y')=g(x',y')\leq g(x^*, y^*)+\lambda \leq g(x,y)+\lambda\leq g(x,y)+\lambda j=f(x,y).
$$
We have shown \eqref{G-min} for the case when $x'=0$.

Step two: we consider the case when $x'\in A_k$, for $k\in \mathbb{N}_\ell$. For all $j\in \mathbb{Z}_{\ell}$,  the second inequality of \eqref{conditions-on-lambda} leads to
$$
g(x',y')+\lambda \ell\leq g(x_j,y_j)+ \lambda j.
$$
This together with the hypothesis that $(x_j,y_j)\in \Omega_j\times\mathbb{R}^{d'}$ is a minimizer of $g$ on $\Omega_j\times\mathbb{R}^{d'}$, for all $j\in \mathbb{Z}_\ell$, ensures that for all $x\in A_j$, $j\in \mathbb{Z}_\ell$ and for all $y\in \mathbb{R}^{d'}$,
$$
f(x',y')=g(x',y')+\lambda k\leq g(x',y')+\lambda \ell\leq g(x_j,y_j)+\lambda j\leq g(x,y)+\lambda j=f(x,y).
$$
For all $x\in A_\ell$, since $(x',y')$ is a minimizer of $g$ on $\Omega_\ell\times\mathbb{R}^{d'}$ and $k\leq \ell$, we find that for all $x\in A_\ell$ and for all $y\in \mathbb{R}^{d'}$,
$$
f(x',y')=g(x',y')+\lambda k\leq g(x,y)+\lambda k\leq g(x,y)+\lambda \ell=f(x,y).
$$
For all $x\in A_j$, for $j=\ell+1, \ell+2, \dots, d$, by using the first inequality of \eqref{conditions-on-lambda} and again the fact that $(x^*,y^*)$ is a global minimizer of $g$, we derive that for all $x\in A_j$, for $j=\ell+1, \ell+2, \dots, d$, and for all $y\in\mathbb{R}^{d'}$,
$$
f(x',y')=g(x',y')+\lambda k\leq g(x^*,y^*)+\lambda(k+1)\leq g(x,y)+\lambda j=f(x,y).
$$
We have shown \eqref{G-min} for the cases $x'\in A_k$, for $k\in\mathbb{N}_\ell$.

Summarizing the above two steps of verification, we conclude that \eqref{G-min} holds true and thus, $(x',y')$ is a global minimizer of $f$ on $\mathbb{R}^d\times\mathbb{R}^{d'}$.
\end{proof}

Similarly to the comment made after the proof of Theorem \ref{choice-of-lambda2}, we now remark on condition \eqref{conditions-on-min}. It is straightforward to verify that condition \eqref{conditions-on-min} is equivalent to
\begin{equation}\label{conditions-on-min***}
g(x',y')\leq \frac{g(x_j,y_j)+(\ell-j)g(x^*,y^*)}{\ell-j+1}, \ \ \mbox{for all}\ \ j\in \mathbb{Z}_\ell.
\end{equation}
Inequality \eqref{conditions-on-min***} show that $g(x',y')$ should be bounded above by an weighted average of $g(x^*,y^*)$ and $g(x_j,y_j)$. By the definition of $(x^*,y^*)$, $(x',y')$ and $g(x_j, y_j)$, we derive that
$$
g(x^*, y^*)\leq g(x',y')\leq g(x_j,y_j).
$$
This shows that condition \eqref{conditions-on-min} is reasonable.


We illustrate  Theorem \ref{choice-of-lambda} by presenting in the next corollary its special case when $\ell=1$. The corollary gives a choice of the parameter $\lambda$ which guarantees a global minimizer of $f$ with sparsity of level 1 for the variable $x$.

\begin{coro}
Let $(x^*,y^*)\in \mathbb{R}^d\times\mathbb{R}^{d'}$ be a global minimizer of $g$, $(x',y')\in \Omega_1\times\mathbb{R}^{d'}$ be a minimizer of $g$ on $\Omega_1\times\mathbb{R}^{d'}$, and $(0,y_0)\in A_0\times \mathbb{R}^{d'}$ be a minimizer of $g$ on $A_0\times\mathbb{R}^{d'}$. Suppose that
\begin{equation}\label{conditions-on-min-ell=1}
g(x',y')-g(x^*,y^*)\leq g(0,y_0)-g(x',y').
\end{equation}
If the parameter $\lambda$ is chosen to satisfy the condition
\begin{equation}\label{conditions-on-lambda-ell=1}
g(x',y')-g(x^*,y^*)\leq \lambda\leq g(0,y_0)-g(x',y'),
\end{equation}
then $(x',y')$ is a global minimizer of $f$ on $\mathbb{R}^d\times\mathbb{R}^{d'}$.
\end{coro}
\begin{proof}
This result is obtained by specializing Theorem \ref{choice-of-lambda} to the special case where $\ell=1$ and noticing that $A_0=\{0\}$. In this case, we have that $x_0=0$.
\end{proof}

The next theorem connects a global minimizer of $g$ with  a global (or local) minimizer of $f$.

\begin{theorem}
Let $(x^*,y^*)\in \mathbb{R}^d\times\mathbb{R}^{d'}$ be a global minimizer of $g$ on $\mathbb{R}^d\times\mathbb{R}^{d'}$.

(i) If $x^*=0$, then $(x^*,y^*)$ is a global minimizer of $f$ on $\mathbb{R}^d\times\mathbb{R}^{d'}$.

(ii) If for some $\ell\in \mathbb{N}_d$, $x^*\in A_\ell$ and
\begin{equation}\label{lowBoundofg(x,y^*)}
g(x^*, y^*)+\lambda(\ell-j)\leq g(x, y),\ \ \mbox{for all}\ \ (x,y)\in A_{j}\times\mathbb{R}^{d'}, \ j\in\mathbb{N}_\ell,
\end{equation}
then $(x^*,y^*)$ is a global minimizer of $f$ on $\mathbb{R}^d\times\mathbb{R}^{d'}$.

(iii) If $x^*\in A_d$, and for some $j\in \mathbb{Z}_{d+1}$ and for some $x'\in A_j$,
\begin{equation}\label{UpperBoundofg(x',y^*)}
g(x^*, y^*)+\lambda(d-j)>g(x', y^*),
\end{equation}
then $(x^*,y^*)$ is a local minimizer of $f$ but not a global minimizer of $f$.
\end{theorem}
\begin{proof}
For both Items (i) and (ii), we shall show that
\begin{equation}\label{f-min-d}
f(x^*, y^*)\leq f(x, y), \ \ \mbox{for all} \ \ (x,y)\in \mathbb{R}^d\times\mathbb{R}^{d'}.
\end{equation}
Since $(x^*,y^*)\in \mathbb{R}^d\times\mathbb{R}^{d'}$ is a global minimizer of $g$, we have that
\begin{equation}\label{g-min-d}
g(x^*,y^*)\leq g(x,y), \ \ \mbox{for all} \ \ (x,y)\in \mathbb{R}^d\times\mathbb{R}^{d'}.
\end{equation}

(i) If $x^*=0$, by the definition of $\|\cdot\|_0$, we have that $\|x^*\|_0=0$ and thus,
$
f(x^*,y^*)=g(x^*,y^*).
$
Using this equation together with condition \eqref{g-min-d}, we observe for all $(x,y)\in
A_j\times\mathbb{R}^{d'}$, for $j\in \mathbb{Z}_{d+1}$ that
$$
f(x^*,y^*)=g(x^*,y^*)\leq g(x,y)\leq g(x,y)+\lambda j=f(x,y).
$$
This together with \eqref{decompositionOfRd} ensures that \eqref{f-min-d} holds.

(ii) Since for some $\ell\in \mathbb{N}_d$, there holds $x^*\in A_\ell$ and since condition \eqref{lowBoundofg(x,y^*)} holds, we have for all $(x,y)\in A_j\times\mathbb{R}^{d'}$, for $j\in \mathbb{Z}_{\ell}$ that
$$
f(x^*,y^*)=g(x^*,y^*)+\lambda\ell=g(x^*,y^*)+\lambda(\ell-j)+\lambda j\leq g(x,y)+\lambda j=f(x,y).
$$
Moreover, by \eqref{g-min-d} we have for $(x,y)\in
A_{\ell+j}\times\mathbb{R}^{d'}$, $j\in \mathbb{Z}_{d-\ell+1}$ that
$$
f(x^*,y^*)=g(x^*,y^*)+\lambda\ell\leq g(x,y)+\lambda(\ell+j)=f(x,y).
$$
Again, according to \eqref{decompositionOfRd}, we find that \eqref{f-min-d} holds.

(iii) Since $x^*\in A_d$ and $(x^*,y^*)$ is a global minimizer of $g$ on $\mathbb{R}^d\times\mathbb{R}^{d'}$,
we have for all $(x,y)\in A_d\times\mathbb{R}^{d'}$ that
$$
f(x^*,y^*)=g(x^*,y^*)+\lambda d\leq g(x,y)+\lambda d=f(x,y).
$$
That is,
\begin{equation}\label{local}
f(x^*,y^*)\leq f(x,y), \ \ \mbox{for all}\ \ (x,y)\in A_d\times\mathbb{R}^{d'}.
\end{equation}
According to Proposition \ref{sets} (iv), we know that $A_d$ is an open set. Thus, $A_d\times\mathbb{R}^{d'}$ is an open set. Therefore, inequality \eqref{local} ensures that the pair $(x^*,y^*)$ is a local minimizer of $f$.

It remains to show that $(x^*,y^*)$ is not a global minimizer of $f$. By condition \eqref{UpperBoundofg(x',y^*)}, we have for some $j\in \mathbb{Z}_{d+1}$ and for some $x'\in A_j$ that
$$
g(x^*,y^*)+\lambda d>g(x',y^*)+\lambda j.
$$
This together with the fact $x^*\in A_d$ and $x'\in A_j$ ensures that
$$
f(x^*,y^*)=g(x^*,y^*)+\lambda d >g(x',y^*)+\lambda j= f(x',y^*).
$$
This implies that $(x^*,y^*)$ is not a global minimizer of $f$.
\end{proof}

In the next theorem, we provide properties of a global minimizer of $f$.

\begin{theorem}\label{necessary-condition}
Let $(x^*,y^*)\in\mathbb{R}^d\times\mathbb{R}^{d'}$ be a global minimizer of $f$ on $\mathbb{R}^d\times\mathbb{R}^{d'}$.

(i) If for some $\ell\in \mathbb{Z}_{d+1}$, $x^*\in A_\ell$, then $(x^*,y^*)$ is a minimizer of $g$ on $\Omega_\ell\times\mathbb{R}^{d'}$.

(ii) If $(x^*,y^*)\in\mathbb{R}^d\times\mathbb{R}^{d'}$ is not a global minimizer of $g$, then $x^*\in \Omega_{d-1}$.
\end{theorem}
\begin{proof}
(i) Since for some $\ell\in \mathbb{Z}_{d+1}$, $(x^*,y^*)\in A_\ell\times\mathbb{R}^{d'}$ is
a global minimizer of $f$, by \eqref{f-min-d}, we have that
$$
g(x^*,y^*)+\lambda \ell=f(x^*,y^*)\leq f(x,y)=g(x,y)+\lambda j, \ \ \mbox{for all} \ \ (x,y)\in A_j\times\mathbb{R}^{d'}, \ \ j\in \mathbb{Z}_{\ell+1}.
$$
It follows that
$$
g(x^*,y^*)+\lambda(\ell-j)\leq g(x,y), \ \ \mbox{for all} \ \ (x,y)\in A_j\times\mathbb{R}^{d'}, \ \ j\in \mathbb{Z}_{\ell+1}.
$$
Using this inequality and noting that
$$
\lambda (\ell-j)\geq 0, \ \ \mbox{for}\ \ j\in \mathbb{Z}_{\ell+1},
$$
we have that
$$
g(x^*,y^*)\leq g(x^*,y^*)+\lambda(\ell-j)\leq g(x,y), \ \ \mbox{for all} \ \ (x,y)\in A_j\times\mathbb{R}^{d'},
\ \ \ell\in \mathbb{Z}_{\ell+1}.
$$
This ensures that $(x^*,y^*)$ is a minimizer of $g$ on $\Omega_\ell\times\mathbb{R}^{d'}$.

(ii) We prove this assertion by contradiction. Assume to the contrary that $x^*\notin \Omega_{d-1}$. Since $x^*\in\mathbb{R}^d$,
$
\mathbb{R}^d=\Omega_{d-1}\cup A_d
$
and
$
\Omega_{d-1}\cap A_d=\emptyset,
$
we must have that $x^*\in A_d$. By Item (i) of this theorem with $\ell=d$, we conclude that $(x^*,y^*)$ is a minimizer of $g$ on $\Omega_d\times\mathbb{R}^{d'}$. Noting that $\Omega_d=\mathbb{R}^d$, we confirm that $(x^*,y^*)$ is a global minimizer of $g$ on $\mathbb{R}^d\times\mathbb{R}^{d'}$, a contradiction. This contradiction ensures that $x^*\in \Omega_{d-1}$.
\end{proof}

Theorem \ref{necessary-condition} (ii) provides a sufficient condition which guarantees that a global minimizer of $f$ is sparse.
The next corollary follows immediately from Theorem \ref{necessary-condition} (ii).

\begin{coro}
If $(x^*,y^*)\in\mathbb{R}^d\times\mathbb{R}^{d'}$ is a global minimizer of $f$,  then, either
$x^*\in \Omega_{d-1}$ or $(x^*,y^*)$ is a global
minimizer of $g$ on $\mathbb{R}^d\times\mathbb{R}^{d'}$.
\end{coro}



We next present an understanding of the relation between the local minimizers of minimization problem \eqref{opt:f} for a fixed parameter $\lambda>0$ and the constrained minimization problem without the term involving the $\ell_0$-norm.
We now define precisely the constrained minimization problem. 
For a given index set $\mathcal{I}$, we introduce a minimization problem on $\mathcal{C}_{\mathcal{I}}\times\mathbb{R}^{d'}$ by
\begin{equation}\label{nonconvex_modelwithout-l0}
\min\left\{
g(x,y):\ \ (x, y)\in \mathcal{C}_{\mathcal{I}}\times\mathbb{R}^{d'}\right\}.
\end{equation}

We need a technical lemma to compare the support of a given vector with that of vectors in its close neighbourhood.

\begin{lemma}\label{equal-support}
If $x^*\in \mathbb{R}^d$ is given, then there exists $\delta_0>0$ such that

(i) for all $x\in B(x^*, \delta_0)$, there holds $S(x^*)\subseteq S(x)$;

(ii) for all $x\in B(x^*, \delta_0)\cap \mathcal{C}_{\mathcal{I}}$ with $\mathcal{I}:=S(x^*)$, there holds $S(x^*)= S(x)$.
\end{lemma}
\begin{proof}
For a fixed number $0<\mu\leq1/2$, we let $\delta_0:=\min\{\mu|x^*_j|: j\in S(x^*)\}.$ Clearly, $\delta_0>0$. Suppose that $i\in S(x^*)$. For all $x\in B(x^*, \delta_0)$, we have that
$$
|x_i^*-x_i|\leq \|x-x^*\|_2\leq \delta_0
$$
and
$$
|x_i|\geq |x^*_i|-|x^*_i-x_i|\geq \delta_0>0.
$$
This implies that $i\in S(x)$. Thus, $S(x^*)\subseteq S(x)$, which proves Item (i).

To show Item (ii), we note that when $x\in \mathcal{C}_{\mathcal{I}}$, by the definition of  $\mathcal{C}_{\mathcal{I}}$, there holds $S(x)\subseteq S(x^*)$. This together with Item (i) yields  $S(x^*)= S(x)$.
\end{proof}

A pair  $(x^*, y^*)\in \mathbb{R}^d\times\mathbb{R}^{d'}$ is called a local minimizer of the minimization problem \eqref{opt:f}, if there exists a $\delta>0$ such that
$$
f(x^*, y^*)\leq f(x,y),\ \ \mbox{for all }\ \ x\in B(x^*, \delta),\ y\in B(y^*,\delta).
$$
Here comes the theorem concerning the relation between local minimizers of minimization problems \eqref{opt:f} and \eqref{nonconvex_modelwithout-l0}.

\begin{theorem}\label{relation:f-and-g}
Suppose that $(x^*, y^*)\in\mathbb{R}^d\times \mathbb{R}^{d'}$ is given.
The pair  $(x^*, y^*)$ is a local minimizer of the minimization problem \eqref{opt:f} with a fixed parameter $\lambda>0$ if and only if  $(x^*, y^*)$ is a local minimizer of the constrained minmization problem \eqref{nonconvex_modelwithout-l0} with $\mathcal{I}:=S(x^*)$.
\end{theorem}
\begin{proof}
Suppose that the pair $(x^*, y^*)$ is a local minimizer of the minimization problem \eqref{opt:f} with a fixed parameter $\lambda>0$ and we show that the pair is a local minimizer of the constrained minmization problem \eqref{nonconvex_modelwithout-l0} with $\mathcal{I}=S(x^*)$.  We prove this by contradiction.
Since $\mathcal{I}:=S(x^*)$, we note that  $x^*\in \mathcal{C}_{\mathcal{I}}$. Assume to the contrary that
the pair $(x^*, y^*)$ is not a local minimizer of the constrained minimization problem \eqref{nonconvex_modelwithout-l0} with $\mathcal{I}:=S(x^*)$. According to the definition of the local minimizer of $g$, we observe that for any $\delta>0$, there exist $x_\delta\in B(x^*,\delta)\cap \mathcal{C}_{\mathcal{I}}$ and $y_\delta\in B(y^*, \delta)$ such that
$
g(x^*, y^*)>g(x_\delta, y_\delta).
$
Item (ii) of Lemma \ref{equal-support} ensures that there exists $\delta_0>0$ such that for all $x\in B(x^*, \delta_0)\cap \mathcal{C}_{\mathcal{I}}$ with $\mathcal{I}:=S(x^*)$, there holds $S(x^*)= S(x)$. Hence, for any $0<\delta<\delta_0$, there exist $x_\delta\in B(x^*,\delta)\cap \mathcal{C}_{\mathcal{I}}$ and $y_\delta\in B(y^*, \delta)$ such that  $S(x^*)= S(x_\delta)$, which implies $\|x^*\|_0=\|x_\delta\|_0$, and
$g(x^*, y^*)>g(x_\delta, y_\delta)$.
This implies that for the fixed parameter $\lambda>0$, there holds
$$
f(x^*,y^*)=g(x^*, y^*)+\lambda\|x^*\|_0>g(x_\delta,y_\delta)+\lambda \|x_\delta\|_0=f(x_\delta,y_\delta).
$$
This violates the assumption that $(x^*,y^*)$ is a local minimizer of the minimization problem \eqref{opt:f} with the parameter $\lambda>0$.

Now, suppose that $(x^*, y^*)$ is a local minimizer of the constrained minimization problem \eqref{nonconvex_modelwithout-l0} with $\mathcal{I}:=S(x^*)$ and we prove that $(x^*, y^*)$ is a local minimizer of the minimization problem \eqref{opt:f} with the fixed parameter $\lambda>0$. We proceed the proof by considering two cases $x\in \mathcal{C}_{\mathcal{I}}$ and $x\notin \mathcal{C}_{\mathcal{I}}$ separately. We first consider the case when $x\in \mathcal{C}_{\mathcal{I}}$. The definition of the local minimizer ensures that there exists a $\delta_1>0$ such that
\begin{equation}
g(x^*, y^*)\leq g(x,y), \ \ \mbox{for all}\ \ x\in B(x^*, \delta_1)\cap \mathcal{C}_{\mathcal{I}}, \ y\in B(y^*,\delta_1).
\label{minimizer-of-g}
\end{equation}
By Item (ii) of Lemma \ref{equal-support}, we have that there exists a $\delta_0>0$ such that for all $x\in B(x^*, \delta_0)\cap \mathcal{C}_{\mathcal{I}}$,  there holds $S(x^*)= S(x)$. This implies that for all $x\in B(x^*, \delta_0)\cap \mathcal{C}_{\mathcal{I}}$, there holds $\|x^*\|_0=\|x\|_0$.  Choose $\delta:=\min\{\delta_0, \delta_1\}$. Then, for all $x\in B(x^*, \delta)\cap \mathcal{C}_{\mathcal{I}}$ and for all $y\in B(y^*, \delta)$, by \eqref{minimizer-of-g},  there holds
$$
f(x^*, y^*)=g(x^*, y^*) +\lambda \|x^*\|_0\leq g(x,y)+\lambda \|x^*\|_0
=  g(x,y)+\lambda \|x\|_0=f(x,y).
$$
This yields
\begin{equation}
f(x^*, y^*)\leq f(x,y),  \ \ \mbox{for all}\ \ x\in B(x^*, \delta)\cap \mathcal{C}_{\mathcal{I}}, \ y\in B(y^*,\delta).
\label{case1}
\end{equation}

We next consider the case when $x\notin \mathcal{C}_{\mathcal{I}}$. By Item (ii) of Proposition \ref{zeroNormCompare}, we conclude that there exists $\delta_2>0$ such that
\begin{equation}\label{comparethel0norm}
    \|x\|_0\geq \|x^*\|_0+1, \ \ \mbox{for all}\ \ x\in B(x^*, \delta_2)\setminus \mathcal{C}_{\mathcal{I}}.
\end{equation}
Since $g$ is continuous, for $\lambda>0$, there exists $\delta_3>0$ such that
\begin{equation}\label{comparevalueofg}
g(x^*, y^*)\leq g(x,y)+\lambda, \ \ \mbox{for all}\ \ x\in B(x^*, \delta_3)\setminus \mathcal{C}_{\mathcal{I}}, \ y\in B(y^*, \delta_3).
\end{equation}
Choose $\delta:=\min\{\delta_j: j=0,1,2,3\}$.
By employing inequality \eqref{comparevalueofg} and then inequality \eqref{comparethel0norm}, we have for all $x\in B(x^*, \delta)\setminus \mathcal{C}_{\mathcal{I}}$, $y\in B(y^*, \delta)$, that
$$
f(x^*,y^*)=g(x^*,y^*)+\lambda\|x^*\|_0\leq g(x,y)+\lambda\|x^*\|_0+\lambda\leq g(x,y)+\lambda\|x\|_0=f(x,y).
$$
That is,
\begin{equation}
f(x^*, y^*)\leq f(x,y),  \ \ \mbox{for all}\ \ x\in B(x^*, \delta)\setminus \mathcal{C}_{\mathcal{I}}, \ y\in B(y^*,\delta).
\label{case2}
\end{equation}
Combining inequalities \eqref{case1} and  \eqref{case2} leads to
$$
f(x^*, y^*)\leq f(x,y),  \ \ \mbox{for all}\ \ x\in B(x^*, \delta), \ y\in B(y^*,\delta),
$$
which implies that the pair $(x^*,y^*)$ is a local minimizer of the minimization problem \eqref{opt:f} with a fixed parameter $\lambda>0$.
\end{proof}

Theorem \ref{relation:f-and-g} is useful in developing efficient numerical algorithms for solving non-convex minimization problems involved functions $f$ having the form \eqref{def:f2} or \eqref{def:f3} and analyzing convergence of the algorithms, since in such cases Theorem \ref{relation:f-and-g} guarantees that the non-convex minimization problems are reduced to convex minimization problems on certain support sets. We next present two corollaries that specialize Theorem \ref{relation:f-and-g} to functions $f$ having a special form \eqref{def:f2} or \eqref{def:f3}.

\begin{coro}\label{Coro:relation:f-and-g}
Suppose that $(x^*, y^*)\in\mathbb{R}^d\times \mathbb{R}^{d'}$ is given.
The pair  $(x^*, y^*)$ is a local minimizer of the minimization problem \eqref{opt:f} with $f$ being defined by \eqref{def:f2} for a fixed parameter $\lambda>0$ if and only if  $(x^*, y^*)$ is a global minimizer of the convex minimization problem \eqref{nonconvex_modelwithout-l0} with $\mathcal{I}:=S(x^*)$ and $g(x,y):=\phi(y)+\mu\|x-Dy\|_2^2$, where $\phi$ is the convex function and $D$ is $d\times d'$ matrix appearing in \eqref{def:f2}.
\end{coro}

The sufficient condition for a pair $(x^*,y^*)$ to be a local minimizer of the minimization problem \eqref{opt:f} presented in Corollary \ref{Coro:relation:f-and-g} for a special example of convex function $g$ was obtained in Proposition 2.3 of \cite{Zeng18}.

\begin{coro}\label{Coro2:relation:f-and-g}
Suppose that $(x^*, y^*)\in\mathbb{R}^d\times \mathbb{R}^{d'}$ is given.
The pair  $(x^*, y^*)$ is a local minimizer of the minimization problem \eqref{opt:f} with $f$ being defined by \eqref{def:f3} for a fixed parameter $\lambda>0$ if and only if  $(x^*, y^*)$ is a global minimizer of the convex minimization problem \eqref{nonconvex_modelwithout-l0} with $\mathcal{I}:=S(x^*)$ and $g(x,y):=\phi(y)+\mu\|x-Dy\|_1$, where $\phi$ is the convex function and $D$ is $d\times d'$ matrix appearing in \eqref{def:f3}.
\end{coro}

\section{Final Remarks}

We briefly discuss possible extension of the results presented in previous sections involving matrix $M$ and comment on potential uses of the main results of this paper.

We first elaborate an extension of the results involving matrix $M$ which has been assumed to satisfy hypothesise  \eqref{MR}. We now suppose that the matrix $M$ has an arbitrary rank $r$ with $0<r\leq \min\{d,m\}$. In this general case, the singular value decomposition of $M$ can be used to remove hypothesise  \eqref{MR} on $M$. Clearly, matrix $M$ has the singular value decomposition
\begin{equation}\label{SVD}
    M=U\Lambda V^*
\end{equation}
where $U$ is a $d\times d$ unitary matrix, $V$ is an $m\times m$ unitary matrix and $\Lambda$ is a $d\times m$ diagonal matrix having the nonzero diagonal entries $\lambda_1\geq \cdots \geq \lambda_r>0$. We can first extend the results in sections 3 and 4 involving matrix $M$ to the case when $M=\Lambda$. The regularization problem \eqref{regularization} with $f$ having the form \eqref{main-function-transform} for this special case is to impose the regularization only for the first $r$ components of the variable $x$ and leave its remaining $m-r$ components not regularized. 
Results for the general case can be obtained by using appropriate changes of variables with the two unitary matrices $U$ and $V$ from the singular value decomposition \eqref{SVD}. We would leave details of the extension to the interested reader.

We have proved rigorously that if the regularization parameter is chosen appropriately, the $\ell_0$ norm regularization will lead to sparse solutions, a result previously validated empirically. 
Regularization parameter choice strategies presented in sections 3 and 4 are all theoretical since in general it is not realistic to know a global minimizer of function $g$. Nevertheless, these results provide insights into the connection between the choice of the regularization parameter with the locations of global minimizers of $f$. They can serve as a guidance for further designing practical parameter choice strategies. For example, one may estimate a global minimizer of $g$ via certain means. In such a case, our parameter strategies may lead to practical uses. This requires further investigation. 

Finally, we indicate that Theorem \ref{relation:f-and-g} and especially Corollaries \ref{Coro:relation:f-and-g} and \ref{Coro2:relation:f-and-g} are useful in developing efficient algorithms for finding local minimizers of the regularized non-convex optimization problems. The essence of Corollaries \ref{Coro:relation:f-and-g} and \ref{Coro2:relation:f-and-g} is that they identify a local minimizer of a non-convex optimization problem with that of a convex optimization problem. Hence, finding a local minimizer of a non-convex optimization problem can be done by finding a local minimizer of a convex optimization problem. In general, solving a convex optimization problem is much easier than solving a non-convex optimization problem.

\end{document}